\documentclass{article}

\usepackage{arxiv}
\usepackage[utf8]{inputenc} 
\usepackage[T1]{fontenc}    
\usepackage{hyperref}       
\usepackage{url}            
\usepackage{booktabs}       
\usepackage{amsfonts}       
\usepackage{nicefrac}       
\usepackage{microtype}      
\usepackage{lipsum}
\usepackage{amsthm,amsmath}
\usepackage{indentfirst}
\usepackage{graphicx}
\usepackage{epstopdf}
\usepackage{fourier}
\usepackage{bm}
\usepackage{mathtools}
\usepackage{latexsym,enumerate}
\usepackage{xcolor}

\def\bm{\boldsymbol}

\newcommand{\BEA}{\begin{eqnarray}}
\newcommand{\EEA}{\end{eqnarray}}

\newcommand{\peq}{p_{eq}}

\newcommand{\td}{\operatorname{d}}

\newcommand{\mbr}{\bm{r}}
\newcommand{\mbd}{\bm{d}}
\newcommand{\mbv}{\bm{v}}
\newcommand{\mbp}{\bm{p}}
\newcommand{\mbq}{\bm{q}}

\newcommand{\opre}{\operatorname{Re}}
\newcommand{\opim}{\operatorname{Im}}
\DeclareMathOperator*{\argmin}{arg\,min}
\DeclareMathOperator*{\diag}{diag}

\newtheorem{thm}{Theorem}[section]

\newtheorem{defn}[thm]{Definition}

\newtheorem{rem}[thm]{Remark}

\title{Linear Response Based Parameter Estimation in the Presence of Model Error}

\author{
  He Zhang \\
  Department of Mathematics \\
  The Pennsylvania State University, University Park, PA 16802, USA.\\
  \texttt{hqz5159@psu.edu} \\
  \And
  John Harlim \\
  Department of Mathematics, Department of Meteorology and Atmospheric Science, Institute for CyberScience\\
  The Pennsylvania State University, University Park, PA 16802, USA.\\
  \texttt{jharlim@psu.edu} \\
  \And
  Xiantao Li \\
  Department of Mathematics\\
  The Pennsylvania State University, University Park, PA 16802, USA.\\
  \texttt{xxl12@psu.edu} \\
}

\begin{document}

\maketitle

\begin{abstract}
Recently, we proposed a method to estimate parameters of stochastic dynamics based on the linear response statistics. The method rests upon a nonlinear least-squares problem that takes into account the response properties that stem from the Fluctuation-Dissipation Theory. In this article, we address an important issue that arises in the presence of model error. In particular, when the equilibrium density function is high dimensional and non-Gaussian, and in some cases, is unknown, the linear response statistics are inaccessible. We show that this issue can be resolved by fitting the imperfect model to appropriate \emph{marginal linear response statistics} that can be approximated using the available data and parametric or nonparametric models. The effectiveness of the parameter estimation approach is demonstrated in the context of molecular dynamical models (Langevin dynamics) with a non-uniform temperature profile, where the modeling error is due to coarse-graining, and a PDE (non-Langevin dynamics) that exhibits spatiotemporal chaos, where the model error arises from a severe spectral truncation. In these examples, we show how the imperfect models, the Langevin equation with parameters estimated using the proposed scheme, can predict the nonlinear response statistics of the underlying dynamics under \emph{admissible external disturbances}.
\end{abstract}


\section{Introduction} \label{sec:intro}

An important and routine task in scientific research is to determine parameters in a given mathematical model. Whether the model is formulated from direct empirical observations or deduced from an underlying dynamic that is much more comprehensive, model error is inevitable in general. Even in the context of first-principle approaches, e.g., quantum-mechanical descriptions \cite{Kohn1965,runge1984density}, they are rarely implemented in its full form. Instead, various reductions are introduced so that the model fits into practical computations, e.g., tight-binding methods \cite{haydock1972electronic} and pseudo-potentials \cite{hamann1979norm}, which necessarily introduce model error,  and the parameters have to be determined in the presence of the model error. 

We have in mind a class of stochastic models, whose parameters can be classified as follows. We assume that there is a subspace of parameters directly responsible for the equilibrium statistics. This subspace of parameters can be determined by many standard statistical methods. Once the parameters in this subspace are determined, the full model parameters can be estimated by fitting the parameters in the quotient space corresponding to the dynamics' response properties. In this paper, we consider fitting the dynamical behavior through an impulse/response approach by perturbing the system with a small external forcing. This approach is known as the linear response in statistical physics as a first step to study non-equilibrium phenomena \cite{Toda-Kubo-2}. A hallmark in linear response theory is the Fluctuation-Dissipation Theory (FDT). Basically, it states that for a dynamical system at equilibrium, the full response to small external perturbations can be characterized by first-order linear responses through appropriate time correlation functions of the {\it unperturbed} dynamical system. Motivated by the accessibility of FDT, we have developed a parameter estimation approach \cite{HLZ:17} as well as an efficient numerical scheme \cite{HLZ:19} under the perfect model assumption. By inferring the estimates from the linear response statistics, we are able to reproduce both the equilibrium statistics and the response properties under small perturbations.

This paper focuses primarily on the response-based method in the presence of model error. In general, the model error may be a result of an empirical assumption, truncation, coarse-graining, multiscale expansion, discounting the memory effect, etc. In practice, the highly desirable parameter estimation methods are those that are not aware of these effects. Such a property is particularly relevant when the only available information is a data set of some observables of interest, which poses several new challenges for parameter estimation methods. First, the external forcing, when applied to the underlying dynamics and then carried over to the imperfect model, may depend on variables other than the observed components. We clarify this issue by introducing a concept called \emph{admissible external forcing}. Another emerging issue is that the linear response involves a variable conjugate to the external force \cite{Toda-Kubo-2} (also interpreted as a flux variable \cite{evans2008statistical}), which, in general, is only accessible when the statistics of the full model are known. We circumvent this difficulty by introducing the marginal linear response statistics that are practically computable given the available data and parametric/nonparametric modeling.

These ideas are best explained with concrete examples. Our first example is motivated by an important problem in computational chemistry,  i.e., coarse-graining (CG) molecular dynamics (MD). The need for CGMD is driven by the observation that full molecular models are limited by the small intrinsic time scale and it is difficult for direct simulations to reach the time scale of interest \cite{baaden_coarse-grain_2013,leach2001molecular,voth2008coarse}, where important biological processes occur \cite{Schlick2002}. CGMD circumvents the problem by projecting the dynamics to the CG variables, a much smaller set, often defined by the local averages of the atomic coordinates and momenta \cite{zhang2008systematic,Noid:13,Noid:08}. The interactions among the CG variables are represented via the free energy, which is also known as the potential of mean force (PMF). In the uniform temperature case, this characterization provides an explicit parametric form of the probability density function. Then traditional statistical methods can be used to determine the parameters in the PMF \cite{chipot2007free}, which we refer to as the equilibrium parameters in the remainder of this article. In addition to this straightforward setup, we will also consider a more general case where the system temperature is non-uniform, which is more challenging since the functional form of the equilibrium density is unknown.

Meanwhile, there are parameters in the CG dynamics that do not appear in the PMF, a well-known example of which is the damping coefficient in the Langevin equation, and more generally, the friction kernel in the generalized Langevin equations. In the latter case, parameters have to be determined by following the statistical properties of the trajectories of the CG variables. Various techniques have been considered, including the Bayesian approach \cite{dequidt2015bayesian}, relative entropy \cite{katsoulakis2013information}, as well as the Pad\'e-type of algorithms that match the time correlation properties \cite{LBL:16,MLL:19,MLC:16}. The application of our response-based approach for CGMD models has been motivated by the steered MD simulations \cite{best2001can,lu1999steered}, where forces are applied to proteins to facilitate the folding and unfolding processes. With a one-dimensional chain model as a simple test problem, we demonstrate how the response-based approach is formulated and implemented to determine the damping coefficient, which is assumed to be a more general band matrix.

Another important class of problems is PDEs with solutions that exhibit chaotic behavior. One well-known example is turbulence. These dynamics are often projected to Fourier modes to study the energy transfer among different modes. From a reduced-order modeling viewpoint, a challenging problem is to deduce a closure model that retains the statistical properties of a finite number of modes. Part of the challenge stems from the fact that the full model's equilibrium statistic is high-dimensional and unknown. We will consider the Kuramoto-Sivashinsky (KS) equation as an example and demonstrate how to circumvent this difficulty with a semi-parametric method using the kernel embedding method. 

In principle, the linear response theory is a prediction tool by itself, and the imperfect models would seem to be of little value in this case.  However, as we previously alluded to, the conjugate variables might not be computable. The semi-parametric approach addresses precisely this issue. Besides, we demonstrate with examples that our approach yields models that also predict the {\it full response} with reasonable accuracy.

The paper is organized as follows: We start with a short review of the linear response theory and the parameter estimation approach based on the linear response in Section \ref{sec:review}. To motivate a general parameter estimation scheme in the presence of model error, a linear fast-slow model will be presented in Section \ref{sec:insight} for which explicit calculations can be carried out. In Section \ref{sec:scheme}, we elaborate on the general difficulties and outline a parameter estimation scheme to address those issues. The numerical scheme will be followed by two examples: a one-dimensional chain model with spatially varying temperature (Section \ref{sec:chain}) and the KS equation (Section \ref{sec:KS}). We close the paper with a summary and discussions in Section \ref{sec:sum}. We include three appendices, discussing the detailed computational method for verifying the results, the derivation of the parametric form of the equilibrium density in the CGMD example, and presenting additional numerical results of the isothermal CGMD.

\section{The parameter estimation method for the perfect model via linear response statistics} \label{sec:review}

Fluctuation-Dissipation Theory (FDT) is a mathematical framework for quantifying the linear response of a dynamical system subjected to small external forcing \cite{leith1975climate}. The linear response statistics, determined based on two-point equilibrium statistics of the unperturbed dynamics, provide estimates for the response properties. In statistical mechanics literature, FDT is known as a linear response approach \cite{evans2008statistical}, which serves as a foundation for defining transport coefficients, e.g., viscosity, diffusion constant, heat conductivity, etc.

We begin by reviewing the linear response theory and the concept of the essential statistics, which is a core component of the response-based estimation method introduced in \cite{HLZ:17}. This will be presented in the context of a stochastic dynamics, expressed in the form of an $n$-dimensional (time-homogeneous) stochastic differential equations (SDEs), also known as the It\^o diffusion \cite{Pav_book:14}. The SDEs, together with its perturbed dynamics, are written as follows,
\BEA
\dot X &=& b(X;\theta) + \sigma(X;\theta)\dot W_t,  \label{eq:Ito}\\
\dot X^{\delta} &=& b(X^{\delta}; \theta) +  c(X^{\delta}) \delta f(t) + \sigma(X^{\delta}; \theta) \dot U_{t}, \label{eq:Ito_per}
\EEA
respectively, where $W_{t}$ and $U_t$ are standard Wiener processes. In the unperturbed system \eqref{eq:Ito}, the vector field $b(X;\theta)$ denotes the drift and $\sigma(X; \theta)$ is the diffusion tensor; while in the perturbed system \eqref{eq:Ito_per}, an order-$\delta$ ($ |\delta|\ll 1$) external perturbation is introduced in the form of $f(x,t) = c(x)\delta f(t)$. In both equations, $\theta\in D$ denotes the model parameters with parameter domain $D\subset \mathbb{R}^{N}$. 

Throughout this section, which concerns with the \emph{perfect model} case, our goal is to estimate the underlying true parameter value $\theta^{\dagger}$ based on the time series of the unperturbed dynamics \eqref{eq:Ito} at the equilibrium. In the case of \emph{imperfect model}, which is the central problem concerned in this paper, there are no true parameter values since the parameters in the imperfect model are not consistent with the parameters in the true model. Before we introduce the concept of the response statistics, we pose the following assumptions for all $\theta\in D$:
\begin{enumerate}
	\item The system governed by Eq. \eqref{eq:Ito} is ergodic with equilibrium density $p_{eq}(x;\theta)$, that is, $p_{eq}(x;\theta)$ is the unique steady-state solution of the corresponding Fokker-Planck equation. In particular, we denote $p_{eq}^{\dagger}(x):= p_{eq}(x;\theta^{\dagger})$. 
	\item The statistical properties associated with Eq. \eqref{eq:Ito_per} can be characterized by a perturbed density $p^{\delta}(x,t;\theta)$, which follows the corresponding Fokker-Planck equation  under initial condition $p^{\delta}(x,0;\theta) = p_{eq}(x;\theta)$.
\end{enumerate}
There is an abuse of notation since for certain SDEs, $p_{eq}$ can depend on a subspace of the parameter. We use the notation $\mathbb{E}_{p_{eq}(\theta)}[\cdot]$ and $\mathbb{E}_{p^\delta(\theta)}[\cdot]$ to emphasis that the corresponding statistics are computed from the sample generated by the corresponding model at parameter value $\theta$ initially at equilibrium.

For any integrable observables $A(\cdot)$, one can define the difference
\begin{equation}\label{eq:full_resp}
\Delta \mathbb{E}[A](t): = \mathbb{E}_{p^{\delta}(\theta)}[A(X^{\delta})](t) - \mathbb{E}_{p_{eq}(\theta)}[A(X)]
\end{equation}
as the full response statistics. The linear response theory allows us to estimate the order-$\delta$ term in \eqref{eq:full_resp} by a convolution integral, that is,
\begin{equation}\label{eq:lin_resp}
\Delta \mathbb{E}[A](t) =\int^{t}_{0}k_{A}(t-s;\theta)\delta f(s) \td s + \mathcal{O}(\delta^2).
\end{equation}
The FDT formulates the linear response operator, $k_{A}(t)$ in \eqref{eq:lin_resp}, as the following two-point statistics:
\begin{equation}\label{eq:lin_oper}
k_{A}(t;\theta): = \mathbb{E}_{p_{eq}(\theta)} \left [A(X(t)) \otimes B \left(X(0); \theta \right) \right], \quad B_{i}(X; \theta):= -\frac{\partial_{X_i}\left[c_{i}(X)p_{eq}(X;\theta)\right]}{p_{eq}(X;\theta)},
\end{equation}
where $B_{i}$ and $c_i$ denote the $i^{\text{th}}$ components of $B$ and $c$, respectively. The variable $B$ will be called the conjugate variable \cite{evans2008statistical} to $A$. We will denote $B(X;\theta^\dagger)$ as $B^\dagger(X)$.

It is worthwhile to point out that in statistical mechanics \cite{zwanzig2001nonequilibrium}, the notion of equilibrium often refers to a statistical ensemble, e.g., the canonical distribution that describes a system at thermal equilibrium. But the approach of linear response can be naturally extended to problems modeled by stochastic differential equations \cite{hairer2010simple}, in which case the equilibrium corresponds to the stationary distribution. The framework of linear response is also applicable to systems at non-equilibrium steady-state \cite{baiesi2009fluctuations}.
 
The significance of FDT is that the response operator is defined without involving the perturbed density $p^{\delta}(x,t;\theta)$. Rather, it can be evaluated at the equilibrium of the unperturbed dynamics. To turn this into a practical tool, we first notice that for a given $t \ge 0$, the value of $k_{A}(t;\theta^{\dagger})$ can be computed using a Monte-Carlo sum based on the time series of the unperturbed system \eqref{eq:Ito} at $p^{\dagger}_{eq}$. For example, let $\left\{X_{i}=X(t_{i})\right\}_{i=1}^{M}$ be the time series generated at $p^{\dagger}_{eq}$ with step length $\Delta t = t_{i+1}-t_{i}$, then for $t = k \Delta t$, the Monte-Carlo approximation can be written as 
\begin{equation}\label{eq:MC_approx}
k_{A}(t; \theta^{\dagger}) \approx \frac{1}{M-k} \sum_{i = 1}^{M-k} A(X_{i+k})\otimes B^{\dagger}(X_i).
\end{equation}
In practice, the right-hand side of \eqref{eq:MC_approx} can be computed directly from data. Meanwhile, the same response statistics can be computed for any $\theta \in D$ by solving \eqref{eq:Ito} over the corresponding parameter value. In particular, 
\begin{equation*}
\hat{k}_{A}(t;\theta):= \mathbb{E}_{p_{eq}(\theta)}[A(X(t)) \otimes B^{\dagger}(X(0))],
\end{equation*}
can still be approximated by a Monte-Carlo sum \eqref{eq:MC_approx} using the time series generated at $p_{eq}(x;\theta)$.

Motivated by the accessibility of $k_{A}(t;\theta^{\dagger})$ as well as the finite-dimensional Pad\'{e} approximation introduced in \cite{LBL:16}, we have proposed to infer the true parameter value $\theta^{\dagger}$ from a discretization of $k_{A}(t;\theta^{\dagger})$. In particular, the discrepancy of the linear response statistics is given by:
\begin{equation*}
\left|\int_{0}^{t}\left(k_{A}^{\dagger}(t-s) - \hat{k}_{A}(t-s; \theta)\right)\delta f(s) \td s  \right| \leq C \left( \int_{0}^{t}\left(k_{A}^{\dagger}(s) - \hat{k}_{A}(s; \theta)\right)^2 \td s  \right)^{\frac{1}{2}},
\end{equation*}
where we have used Cauchy-Schwartz inequality with constant $C$ to denote the $L^2$-norm of the term $\delta f(t)$ in the external forcing.  Minimizing the $L^{2}$-error in the right-hand side on discrete time \emph{essential statistics} $\left\{k^{\dagger}_{A}(t_i):=k_{A}(t_{i};\theta^{\dagger})\right\}_{i=1}^{K}$, 
leads to the nonlinear least-squares problem \cite{HLZ:17,HLZ:19}:
\begin{equation}\label{eq:nonlin_ls}
\hat \theta: = \argmin_{\theta\in D} \sum_{i = 1}^{K} f^{2}_{i}(\theta), \quad f_{i}(\theta): = k^{\dagger}_{A}(t_i) - \hat{k}_{A}(t_{i};\theta), \quad i = 1,2,\dots, K.
\end{equation}
To ensure the solvability of \eqref{eq:nonlin_ls}, we assume that the total number of involved essential statistics is always strictly greater than the dimension of $\theta$, that is, $K>N$. In terms of choosing $\{t_i\}$, one can interpret the optimization problem \eqref{eq:nonlin_ls} as a curve-fitting problem. Thus, the points $\{t_i\}$'s should be selected such that $\{k^{\dagger}_{A}(t_i)\}$'s  represent the qualitative feature in $k^{\dagger}_{A}(t)$. For example, one may include more points around local extreme points of $k^{\dagger}_{A}(t)$. Also, we do not want $t_i$ to be too large, since in that case the value of $k^{\dagger}_{A}(t)$ would be too small and the Monte-Carlo error introduced by \eqref{eq:MC_approx} will overwhelm the significance of the essential statistics.

\begin{rem}
The accessibility of $k_{A}(t;\theta^{\dagger})$ also allows us to use its time derivatives as the essential statistics to infer $\theta^{\dagger}$. For example, we have shown that in \cite{HLZ:17}, under specific external forcing and observables $A(\cdot)$, the first-order time derivative of $k_{A}$ at $t=0$ of a Langevin dynamics model provides direct estimates for the damping coefficients. It is worthwhile to mention that estimating higher-order derivatives of  $k_{A}(t;\theta^{\dagger})$ requires time series with sufficient accuracy, which is not necessarily available, especially when the data are from experimental observations.
\end{rem}

\begin{rem}
It is important to point out that, in general, the quantity $B(X;\theta)$ in \eqref{eq:lin_oper} depends on the unknown parameter $\theta$, which provides difficulties in computing $k_{A}(t;\theta^{\dagger})$ given only a time series of $X$ at $p^{\dagger}_{eq}$. In \cite{Di:16a,Di:16b,Di:18}, Qi and Majda addressed this issue by computing the ``kicked response'' based on a direct simulation of the underlying dynamics. However, since the underlying dynamics is unknown in our scenario and estimating $p^\dagger_{eq}$ is a difficult computational task in high-dimensional problems, we will consider estimating a marginal invariant density $p_R$ defined on the resolved (and identifiable unresolved) variables and define $B$ using the estimated marginal density $p_R$ in placed of $p^{\dagger}_{eq}$, as explained in Section~\ref{sec:scheme}.
\end{rem}

One practical challenge in solving \eqref{eq:nonlin_ls} is that, except for the very special cases, there is no explicit expressions for $\hat{k}_{A}(t;\theta)$. Therefore an iterative method, e.g., the Gauss-Newton method, would necessarily require solving the model \eqref{eq:Ito} repeatedly in a sequential manner, which can be rather computationally demanding.  We have mitigated such an issue by a polynomial based surrogate model in \cite{HLZ:19}, where the cost functions $f_{i}(\theta)$ in \eqref{eq:nonlin_ls} are approximated by linear combinations of orthogonal polynomials, denoted as $f_{i}^{M}(\theta)$, based on the precomputed values $ \left\{f_{i}(\theta_{j})\right\}$ over a set of grid points $\{\theta_{j}\} \subset D$, and such pre-computation can be proceed in parallel. Here, $M$ stands for the order of the polynomials involved. In the subsequent approximations, we will replace $f_{i}(\theta)$ in \eqref{eq:nonlin_ls} by $f_{i}^{M}(\theta)$, and formulate a new least-squares problem for the order-$M$ polynomial based surrogate model. In \cite{HLZ:19}, we have deduced the convergence analysis of such approximation. In particular, we proved that the minimizers $\{\theta^{*}_{M}\}$ of the approximated least-squares problem converge to the minimizer of the true least-squares problem as $M \rightarrow \infty$ under the perfect model assumption and appropriate regularity assumption.

To this end, we have reviewed a parameter estimation approach for It\^o diffusions \eqref{eq:Ito} as well as an efficient numerical scheme under the perfect model assumption. By inferring from the linear response statistics, our estimates can reproduce both the equilibrium statistics and the response under certain perturbations. A much more difficult problem, however, is when the available model is a truncated dynamics or only partially known. In this case, several interesting issues will emerge, which will be investigated in the remaining part of the paper. In the next section, we start with a simple two-scale system to demonstrate the importance of the response statistics in determining the model parameters.

\section{A linear fast-slow model for insight} \label{sec:insight}

Consider the following two-dimensional linear fast-slow dynamics as the underlying model
\BEA
\dot x &=& (a_{11}x + a_{12}y)  + \sigma_{x} \dot W_{x},  \label{eq:slow} \\ 
\dot y &=& \frac{1}{\epsilon}(a_{21}x + a_{22}y)  + \frac{\sigma_{y}}{\sqrt{\epsilon}}\; \dot W_{y}, \label{eq:fast}
\EEA
where $x$ and $y$ are the slow and the fast variables, respectively. Here, $W_x$ and $W_y$ are two independent Wiener processes and $\epsilon \ll 1$ characterizes the time-scale gap. Such multi-scale Ornstein-Uhlenbeck process has been used as a linear test model in a variety of situations \cite{GH:13}. For the sake of non-degeneracy and ergodicity, we assume that
\begin{equation*}
\sigma_x, \sigma_y > 0, \quad a_{11}a_{22} - a_{12}a_{21}>0, \quad a_{11}+ \frac{a_{22}}{\epsilon} <0, \quad a_{22}<0,
\end{equation*}
which are enough to ensure the drift coefficients 
\begin{equation*}
A_{\epsilon}: = \begin{pmatrix}
a_{11} & a_{12} \\
\nicefrac{a_{21}}{\epsilon} & \nicefrac{a_{22}}{\epsilon}  \end{pmatrix}
\end{equation*}
to be negative definite as well as the existence of the averaged (effective) dynamics,
\begin{equation}\label{eq:lin_ave}
\dot X = \tilde{a} X + \sigma_{x} \dot W_t, \quad \tilde{a} := a_{11} - \frac{a_{12}a_{21}}{a_{22}}<0.
\end{equation}
Here, the solution, $X(t)$, with initial condition $X(0)=x(0)$, converges strongly up to order-one time $t$ to the solution $x(t)$ of the true dynamics in \eqref{eq:slow}-\eqref{eq:fast}. See for example \cite{Pav_book:08} for a detailed discussion.

We will keep the form of \eqref{eq:lin_ave} by posing
\begin{equation}\label{eq:lin_red}
\dot X = b X +  \sigma \dot W_t, \quad b<0, \quad \sigma >0,
\end{equation}
as our imperfect model of the fast-slow dynamics with $\theta = (b, \sigma)$ as the unknown parameters to be determined. Applying the parameter estimation scheme reviewed in Section \ref{sec:review} to estimate $\theta$ by taking observable $A(X) = X$ and a time-dependent external forcing $\delta f(t)$, we find that (see the proof of Theorem 1 in \cite{HLZ:17} for computational details),
\begin{equation*}
\hat{k}_{A}(t; \theta) = \mathbb{E}_{p_{eq}(\theta)}[X(t)X(0)]/S^{\dagger} = e^{tb} \frac{S}{S^{\dagger}},
\end{equation*}
where $S = -\sigma^2/2b$ is the equilibrium variance of $X$. In contrast, $S^{\dagger}$ is the variance that corresponds to the ``true parameter'' value $\theta^{\dagger}$, when the underlying dynamics is exactly of the form \eqref{eq:lin_red} with $\theta=\theta^{\dagger}$. Fitting the essential statistics ($k^{\dagger}_{A}(t) = \hat{k}_{A}(t;\theta)$) leads to
\begin{equation}\label{eq:fit1}
\mathbb{E}_{p^{\dagger}_{eq}}[X(t)X(0)]=e^{tb^{\dagger}}S^{\dagger} = e^{tb}S.
\end{equation}
In sharp contrast to the perfect model situation, in our case, the partial observations are from an underlying two-dimensional dynamics \eqref{eq:slow}-\eqref{eq:fast}, while the essential statistics in \eqref{eq:fit1} are derived based on a one-dimensional imperfect model \eqref{eq:lin_red}. Thus, the left-hand side of \eqref{eq:fit1}, legitimately, is not available, and will be replaced by the two-point statistics of the slow variable $x$ in \eqref{eq:slow}. Eq. \eqref{eq:fit1} now becomes,
\begin{equation}\label{eq:fit2}
\mathbb{E}_{p^{\dagger}_{eq}}[x(t)x(0)] = [1,0] e^{tA_{\epsilon}}\Sigma [1,0]^{\top} \approx e^{tb}S = - e^{tb}\frac{\sigma^2}{2b},
\end{equation}
where the covariance matrix $\Sigma = (\Sigma_{ij})$ is determined by the  Lyapunov equation \cite{Pav_book:14}
\begin{equation}\label{eq:lyap}
A_{\epsilon}\Sigma + \Sigma A_{\epsilon}^{\top} = - Q, \quad Q = 
\begin{pmatrix} \sigma_{x}^2 & 0 \\ 0 & \nicefrac{\sigma_{y}^2}{\epsilon} \end{pmatrix}.
\end{equation}

One can interpret Eq. \eqref{eq:fit2} as the consistency between the two time auto-covariance functions of the slow variable $x$ and the variable $X$ in the imperfect  model. Obviously, seeking an exact fit for \eqref{eq:fit2} is not feasible since its left-hand side has a total of $6$ degrees of freedom ($A_\epsilon$ and $\Sigma$), which is much greater than that of the right-hand side. However, taking advantage of the multi-scale structure of the underlying dynamics, we can still find estimates for $\theta$ when $\epsilon\ll 1$ such that \eqref{eq:fit2} is fulfilled up to an error of certain order. To see this, we first simplify the matrix exponential in \eqref{eq:fit2} by computing the eigenvalues of $A_{\epsilon}$,
\begin{equation*}
\lambda_{1,2} = \frac{1}{2} \left(a_{11} + \frac{a_{22}}{\epsilon} \pm \sqrt{\Delta} \right), \quad 
\Delta = \left(a_{11}+\frac{a_{22}}{\epsilon} \right)^{2} - \frac{4a_{22}\tilde{a}}{\epsilon},
\end{equation*}
where the $\tilde{a}$ is the same as in Eq. \eqref{eq:lin_ave}. Assuming that the two eigenvalues are real, via a direct asymptotic expansions with respect to $\epsilon$, one can show that
\begin{equation*}
\lambda_{1} = \tilde{a} - \frac{a_{12}a_{21}\tilde{a}}{a_{22}^{2}}\epsilon + O(\epsilon^2), \quad \lambda_{2} = \frac{a_{22}}{\epsilon} + \frac{a_{12}a_{21}}{a_{22}} + \frac{a_{12}a_{21}\tilde{a}}{a_{22}^{2}}\epsilon + O(\epsilon^2).
\end{equation*}
Notice that $a_{22}<0$ due to our assumption, and we are able to truncate the term containing $e^{t\lambda_2}$ in \eqref{eq:fit2} as long as $t$ is large enough. Thus, up to a negligible error, \eqref{eq:fit2} becomes
\begin{equation}\label{eq:fit3}
e^{\lambda_1 t}\frac{\Sigma_{11}(\lambda_2 - a_{11}) - \Sigma_{21}a_{12}}{\lambda_2- \lambda_1} = - e^{tb}\frac{\sigma^2}{2b},
\end{equation}
which can be achieved by taking
\begin{equation}\label{eq:lin_est}
\begin{split}
\hat{b} &= \lambda_1 = \tilde{a} - \frac{a_{12}a_{21}\tilde{a}}{a_{22}^{2}}\epsilon + O(\epsilon^2), \\
\hat\sigma^2  &= -2\lambda_1\frac{\Sigma_{11}(\lambda_2 - a_{11}) - \Sigma_{21}a_{12}}{\lambda_2- \lambda_1}=\sigma_{x}^2 - \frac{2\sigma_{x}^2a_{12}a_{21}- \sigma_{y}^2a_{12}^2}{a_{22}^2}\epsilon + O(\epsilon^2),
\end{split}
\end{equation}
as our estimates. The same estimates have been obtained based on the Riccati equation in \cite{BH:14}. 

Compared to the averaged dynamics \eqref{eq:lin_ave}, our estimate in \eqref{eq:lin_est}, based on matching the linear response statistics, exhausts the capability of the imperfect model \eqref{eq:lin_red} in the sense that the two-point statistics $\mathbb{E}_{p_{eq}}[x(t)x(0)]$, where $t\gg \epsilon$, of the slow variable $x$ can be reproduced up to any order of $\epsilon$ as $\epsilon \rightarrow 0^{+}$.

\begin{figure}[ht]
	\begin{center}
		\includegraphics[width=.75\textwidth]{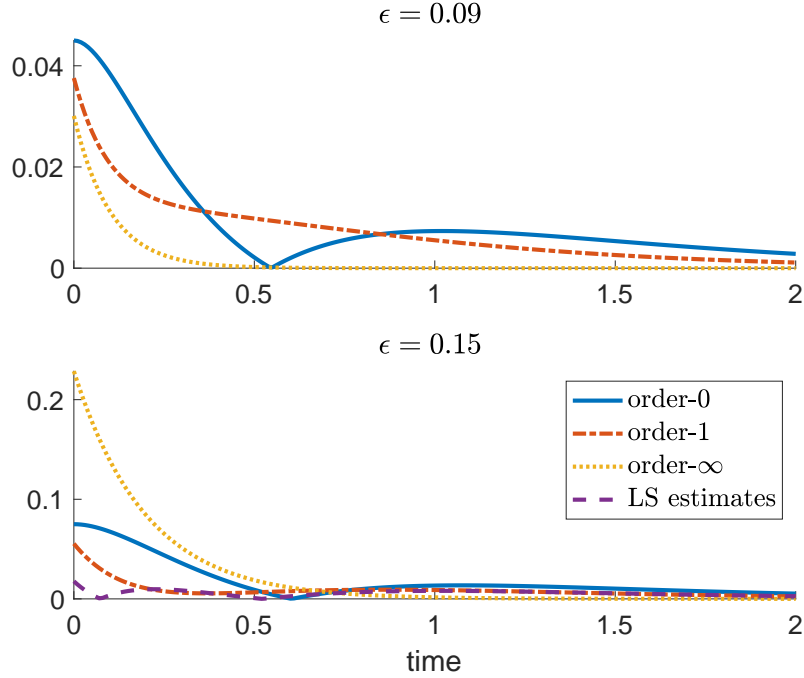}
		\caption{The absolute error in reproducing $\mathbb{E}_{p^{\dagger}_{eq}}[x(t)x(0)]$ under two different scale gaps $\epsilon = 0.09$ (top) and $\epsilon = 0.15$ (bottom). In terms of the estimates, order-$0$: $\hat b = \tilde{a}$ and $\hat \sigma^2=\sigma_{x}^2$; order-1: $\hat{b} = \tilde{a} - \epsilon a_{12}a_{21}\tilde{a}/a_{22}^{2}$ and $\hat{\sigma}^2 = \sigma_{x}^2 - \epsilon(2\sigma_{x}^2a_{12}a_{21}- \sigma_{y}^2a_{12}^2)/a_{22}^2$; order-$\infty$: $\hat{b}$ and $\hat{\sigma}^2$ are given by \eqref{eq:lin_est} without truncation, and LS: least-squares estimates.}
		\label{fig:fast_slow}
	\end{center}
\end{figure}

As a test, we set $a_{11}=a_{21}=a_{22}=-1$, $a_{12}=1$, and $\sigma^2_{x}= \sigma_{y}^2=2$ in \eqref{eq:slow} and \eqref{eq:fast}. Figure \ref{fig:fast_slow} shows the absolute error in reproducing the two-point statistics under two different choices of $\epsilon$, in comparison to various estimation methods of the imperfect model. In particular, we will consider the estimate based on the classical averaging model in \eqref{eq:lin_ave}, which we denote as the order-0 model in Fig.~\ref{fig:fast_slow} since it is the leading order term in \eqref{eq:lin_est}. In addition, we also consider the order-1 model, corresponding to truncating the order-$\epsilon^2$ in \eqref{eq:lin_est}, and the
order-$\infty$ model, corresponding to no truncation in \eqref{eq:lin_est}, i.e., Eq. \eqref{eq:fit3} is fulfilled. For small $\epsilon$, the order-$\infty$ model produces the most accurate estimation. In this regime, we empirically found that the estimation based on parameters obtained by solving Eq.~\eqref{eq:fit2} in the least-squares sense is nearly identical to the order-$\infty$ model (and thus, not reported here).

However, when $\epsilon$ is relatively large, we are no longer allowed to truncate the term $e^{t\lambda_2}$ in the simplification of \eqref{eq:fit2}, which in turn invalidates the expansion in \eqref{eq:lin_est}. One can see the inaccurate estimate of the order-$\infty$ model in the second panel of Figure~\ref{fig:fast_slow}. On the other hand, the solution of \eqref{eq:fit2} in the least-squares sense is accurate and this justifies the importance of solving the nonlinear least-squares problem in \eqref{eq:nonlin_ls} for parameter estimation in the presence of model error.

To summarize, we have elucidated the importance of the proposed parameter estimation scheme in the presence of model error on an idealistic problem. In particular, we have shown that the linear response-based parameter estimation scheme produces accurate linear mean response statistics, outperforming the conventional method, the averaging theory \cite{Pav_book:08}, even in a regime when such theory is valid.

\section{Parameter estimation in the presence of model error} \label{sec:scheme}

Inspired by the result in Section~\ref{sec:insight}, we consider extending the response-based parameter estimation scheme in the presence of model error arising from coarse-graining or missing dynamics. In the current Section, we provide a formal discussion on the key issues in such a difficult scenario and propose a strategy to address them. 

To elaborate the general idea, let us use the abstract notation, 
\begin{equation}\label{eq:full_dyn}
\dot u = F(u), \quad X = P(u),
\end{equation}
to denote the underlying full dynamics with partial observation of the variable of interest, $X$. Problems modeled by SDEs can be treated in a similar manner and they will be addressed later. Here, $P$ (not necessarily available) denotes a projection operator to the observable. We assume that $F$ is unknown, instead, the available imperfect model is given as follows,
\begin{equation}\label{eq:red_dyn}
\dot X = R_{1}(X,Y; \theta), \quad \dot Y = R_{2}(X,Y; \theta).
\end{equation}
Here, we have included an additional variable $Y$ as a possible auxiliary variable and $\theta \in D$ are the unknown model parameters. Our goal is to estimate the parameter $\theta$ in \eqref{eq:red_dyn} using the time series of $X$ observed from \eqref{eq:full_dyn} such that the resulting model in \eqref{eq:red_dyn} can predict the full nonlinear response of the underlying dynamics \eqref{eq:full_dyn} under \emph{admissible external forcing}, as explained below. While the choice of model in \eqref{eq:red_dyn} is central to the accuracy of the prediction, our focus is not on the specification of this model. Instead, we assume that the model in \eqref{eq:red_dyn} is given and our focus, again, is to generalize the proposed parameter estimation scheme reviewed in Section~\ref{sec:review} to this scenario.

In this case, the linear response operator \eqref{eq:lin_oper} associated with the underlying dynamics \eqref{eq:full_dyn} is no longer accessible since, in general, it is a two-point statistic that can  be computed {\it only} if the full data set of $u$ is available. Thus, to infer the parameter $\theta$ in \eqref{eq:red_dyn} from the linear response statistics of the underlying dynamics \eqref{eq:full_dyn}, we need to answer the following two questions:
\begin{enumerate}
\item How to find an alternative to the essential statistics $k_{A}^{\dagger}(t_{i})$ in \eqref{eq:nonlin_ls}, which is not computable given only the time series of $X$?
\item How to formulate a computationally feasible nonlinear least-squares problem analogous to \eqref{eq:nonlin_ls} for the imperfect model \eqref{eq:red_dyn}?
\end{enumerate}

To address these issues, we introduce the concepts of (i) \emph{admissible external forcing}, and (ii) \emph{marginal linear response (MLR)}. Using the notations in \eqref{eq:full_dyn}, similar to \eqref{eq:Ito} and \eqref{eq:Ito_per}, the perturbed underlying dynamics can be written in the form
\begin{equation}\label{eq:full_dyn_per}
\dot u^{\delta} = F(u^{\delta})+ c(u^{\delta})\delta f(t) , \quad X^{\delta} = P(u^{\delta}).
\end{equation}
Taking the time derivative of $X$ in \eqref{eq:full_dyn} and $X^{\delta}$ in \eqref{eq:full_dyn_per}, we have
\begin{equation*}
\dot X = \nabla P(u) F(u), \quad \dot X^{\delta} = \nabla P(u^{\delta}) F(u^{\delta}) + \delta f(t) \nabla P(u^{\delta}) c(u^{\delta}),
\end{equation*}
where $\nabla P$ denotes the Jacobian matrix of $P$. Thus, the term $\delta f(t) \nabla P(u^{\delta}) c(u^{\delta})$ can be interpreted as the external force on the dynamics of $X$. In the case when the full dynamics \eqref{eq:full_dyn} is a system of  SDEs with an external force in the drift term, a similar observation can be made by using the It\^{o}'s formula. To this end, we define:

\begin{defn}\label{def:adm}
An external forcing $f(u,t)= c(u)\delta f(t)$ to the full dynamics in \eqref{eq:full_dyn_per} is admissible with respect to the dynamics $\eqref{eq:full_dyn}$ if  $\;\nabla P(u) c(u)$ can be written as a function of $X$ denoted as $c_{R}(X)$.
\end{defn}

In other words, an admissible external forcing $c(u^{\delta})\delta f(t)$ of the underlying dynamics $\eqref{eq:full_dyn}$ introduces a perturbation $c_{R}(X^{\delta})\delta f(t)$ to the imperfect model \eqref{eq:red_dyn},
\begin{equation}\label{eq:red_dyn_per}
\dot X^{\delta} = R_{1}(X^{\delta},Y^{\delta}; \theta)+ c_{R}(X^{\delta})\delta f(t), \quad \dot Y^{\delta} = R_{2}(X^{\delta},Y^{\delta}; \theta).
\end{equation}
Intuitively, this condition rules out other external forces that depend on $u$, which is not accessible since the imperfect model does not depend explicitly on $u$. We should point out that while the condition in Definition~\ref{def:adm} provides a mathematical justification for the relevant perturbed imperfect model \eqref{eq:red_dyn} that is consistent with the perturbed underlying dynamics in \eqref{eq:full_dyn_per}, this condition is not practically verifiable unless if $P, c$ and $X$ are known. In our numerical examples, we will choose external forcings that are constant and nonzero on resolved components $X$ such that the condition in Definition~\ref{def:adm} is satisfied and they are admissible.

To simplify the discussion, let us define $\theta_{eq}:=\theta_{eq}(\theta)$ as the components of the parameters that directly appear in the equilibrium density $\rho_{eq}(X,Y;\theta_{eq})$ of the imperfect model \eqref{eq:red_dyn}; they usually can be estimated directly by matching equilibrium statistics. For example, in the linear problem in \eqref{eq:lin_red}, while the model parameters $\theta$ consist of the damping coefficient $b$ and noise amplitude $\sigma$, the parameter that appears in the equilibrium density, $\theta_{eq}$, is $S=-\sigma^2/{2b}$, which can be directly obtained by matching the equilibrium variances, i.e., considering $t = 0$ in Eq. \eqref{eq:fit2}.

Once $\theta_{eq}$ has been estimated (by $\hat{\theta}_{eq}$), given an observable $A(X)$, the linear response operator \eqref{eq:lin_oper} of the perturbed imperfect dynamics \eqref{eq:red_dyn_per} is given as,
\begin{equation} \label{eq:lin_oper_red}
\hat{k}_{A}(t;\theta) = \mathbb{E}_{\rho_{eq}(\theta)} \left [A(X(t)) \otimes B\left(X(0),Y(0);\hat\theta_{eq} \right) \right],  \quad B_{i}\left(X,Y;\hat\theta_{eq}\right) = -\frac{\partial_{X_i}\left[\left(c_{R}(X)\right)_{i}\rho_{eq}\left(X,Y;\hat\theta_{eq}\right)\right]}{\rho_{eq}\left(X,Y;\hat\theta_{eq}\right)},
\end{equation}
where the expectation over $\rho_{eq}(\theta)$ will be realized through Monte-Carlo averages with respect to the solutions of the unperturbed imperfect model in \eqref{eq:red_dyn} for parameter $\theta$ that satisfies the constrained of $\theta_{eq}(\theta)=\hat\theta_{eq}$.

Given an observable $A(X)$ and an admissible external forcing $f(u,t)$ in Definition~\ref{def:adm}, in analogous to the linear response operator \eqref{eq:lin_oper} and the linear response statistics \eqref{eq:lin_resp}, we define:
\begin{defn}
Let $p_R(X):=\int \delta(X-P(u)) p^{\dagger}_{eq}(u)\,du$ be the marginal invariant density of $X$, which can be estimated from the given data of $X$ at the equilibrium. We define the marginal linear response (MLR) operator,
\begin{equation}\label{eq:mar_esst}
\hat{k}_{A}^{\dagger}(t):  = \mathbb{E}_{p_{eq}^{\dagger}}\left [A(X(t)) \otimes B^{\dagger}_{R}(X(0)) \right],  \quad \left(B^{\dagger}_{R}\right)_{i}: = -\frac{\partial_{X_i}\left[\left(c_{R}(X)\right)_{i}p_{R}(X)\right]}{p_{R}(X)}.
\end{equation} 
We call $\left\{\hat{k}_{A}^{\dagger}(t_i)\right\}_{i=1}^K$ the marginal essential statistics. Furthermore, we define,
\BEA
\hat\Delta \mathbb{E}[A](t) =\int^{t}_{0}\hat{k}_{A}^\dagger(t-s)\delta f(s) \td s\label{MLR333}
\EEA
as the MLR statistics.
\end{defn}
The two-point statistics in \eqref{eq:mar_esst}, defined with respect to the equilibrium distribution of the underlying dynamics \eqref{eq:full_dyn}, can be approximated via Monte-Carlo as in \eqref{eq:MC_approx}, based on the available time series of $X$. We require the external forcing to be admissible; otherwise, $c_R$ in \eqref{eq:mar_esst} would depend on $u$ and the marginal essential statistics are not accessible when the data of $u$ is not available. In practice, as we will show in the next two sections, the marginal density $p_{R}$ can be learned from the available time series of $X$ by using standard statistical methods. This way, the two-point statistics in \eqref{eq:mar_esst} is computable as opposed to \eqref{eq:MC_approx}. We should also point out that, in principle, $c_R$ can also depend on $Y$. In this case, computing $k_A^\dagger$ in \eqref{eq:mar_esst} requires samples of, both, $X$ and $Y$ for estimating
the marginal density $p_R$, which is a function of $X$ and $Y$. In Section~\ref{sec:reduce_KS}, we will consider the case where the samples data of $Y$ are not available.

To infer the parameter $\theta$ in \eqref{eq:red_dyn} from \eqref{eq:mar_esst}, compared with \eqref{eq:nonlin_ls}, we formulate the following new nonlinear least-squares problem
\begin{equation}\label{eq:new_nls}
\hat \theta : = \argmin_{\left\{\theta\in D\;|\;\theta_{eq}(\theta)=\hat\theta_{eq}\right\}} \sum_{i = 1}^{K} f^{2}_{i}(\theta), \quad f_{i}(\theta): = \hat{k}^{\dagger}_{A}(t_i) - \hat{k}_{A}(t_{i};\theta), \quad i = 1,2,\dots, K,
\end{equation}
where $\hat{k}^{\dagger}_{A}(t_i)$ is the marginal essential statistics \eqref{eq:mar_esst} and $\hat{k}_{A}(t_{i};\theta)$ is the linear response operator \eqref{eq:lin_oper_red} of the imperfect model \eqref{eq:red_dyn}. Eq. \eqref{eq:new_nls} is a dynamic constrained nonlinear least-squares problem through the imperfect model in \eqref{eq:red_dyn}. We should point out the minimization is over a subspace of the parameter set, characterized by the constraint, $\theta_{eq}(\theta)=\hat\theta_{eq}$, where $\hat\theta_{eq}$ is the estimate for the equilibrium parameters, $\theta_{eq}$, obtained from direct matching to the available equilibrium statistics.

In the next two sections, we validate the estimate $\hat{\theta}$ by comparing the nonlinear response of the imperfect model at the estimated parameters with that of the underlying model under admissible external forcing. These statistics are computed via Monte-Carlo averages over solutions of the perturbed underlying dynamics \eqref{eq:full_dyn_per} and perturbed imperfect model \eqref{eq:red_dyn_per}, where the former is only possible when the full underlying dynamic is known. See Appendix~\ref{app:A} for the detailed computational strategy to compute these statistics for stochastic dynamics.  

\begin{rem}\label{rm:1}
For the fast-slow model discussed in Section~\ref{sec:insight}, the underlying dynamics are governed by \eqref{eq:slow}-\eqref{eq:fast}, while the imperfect model is given by \eqref{eq:lin_red}. The time-dependent external forcing we proposed in the previous analysis is admissible since it acts only on the slow variable $x$.
\end{rem}

\section{A one-dimensional molecular model under a local temperature profile}  \label{sec:chain}

We consider a MD model consisting of a chain of atoms with mass $m$. Let the equilibrium spacing between the particles be $a_0$, and the displacement of the $i^{\text{th}}$ particles from its equilibrium position $R_{i}=i a_{0}$ be $r_{i}$. We take the Lennard-Jones (LJ) potential \cite{ishimori:82}, expressed as

\begin{equation}\label{eq:LJ}
U_0(r) = \sum_{i}\sum_{j = i-2}^{i-1} \psi\left(r_i - r_j+(i-j)a_0\right), \quad \psi(r)= \left|r\right|^{-12} - \left|r\right|^{-6}.
\end{equation}

\subsection{The underlying dynamics}
Our underlying dynamics, as a finite-dimensional approximation of the LJ lattice model, contains a total of $N$ particles equipped with periodic boundary conditions. The size of the system is $Na_0$. Let  $\mbr = (r_{1}, r_{2}, \dots, r_{N})^{\top}$ and $\mbv= (v_1, v_2, \dots, v_{N})^{\top}$ be the displacement and the velocities of the particles, respectively. The potential energy of the finite system, as a function of the relative displacement $\mbd: = (r_{2}-r_{1}, r_{3}-r_{2}, \dots, r_{N} - r_{N-1})^{\top}\in \mathbb{R}^{N-1}$, can be written as,
\begin{equation} \label{eq:T_LJ}
U(\mbd)  = \sum_{i=1}^{N+2}\sum_{j = i-2}^{i-1} \psi\left( r_i-r_j+(i-j)a_0\right). 
\end{equation}
The potential $U$ determines the direct interactions among the particles. In the model \eqref{eq:T_LJ}, each particle is interacting with the nearest two neighbors from each side.  We consider a Langevin dynamics of $(\mbr, \mbv)$, driven by the potential energy $U(\mbd)$, the frictional, and random forcing,
\begin{equation}\label{eq:Lan_temp}
\begin{split}
\dot{r}_i &= v_i,\\
\dot{v}_i &= - \nabla_{r_{i}}U(\mbd) - \gamma v_i +\sqrt{2k_{B}T_i \gamma}\dot{W}^{(i)}_{t},
\end{split}
\end{equation}
for $i = 1, 2,\dots, N$, where the mass has been set to unity ($m=1$), $\gamma$ denotes the friction constant, and $W_{t} = \left(W^{(1)}_{t}, \dots, W^{(N)}_t\right)$ is an $N$-dimensional Wiener process. The temperature in \eqref{eq:Lan_temp} is allowed to be non-uniform. As a specific example, we consider a sinusoidal profile
\begin{equation}\label{eq: t-loc}
    \quad k_{B}T_{i} = k_{B}T_{0} + k_{B}\Delta T \sin\left(\frac{2\pi i}{N}  \right), \quad i = 1,2, \dots, N.
\end{equation}
The dynamic in \eqref{eq:Lan_temp} models coupled oscillators interacting with thermal reservoirs of different temperatures (e.g., \cite{lepri2003thermal}). Most of the effort in CGMD models has been focused on systems at a constant temperature ($\Delta T \equiv 0$ in \eqref{eq: t-loc}), in which case, the equilibrium density is known (e.g., \cite{Noid:08}). Unlike the isothermal case, the equilibrium density of \eqref{eq:Lan_temp} is not available in general. We will discuss how this issue is addressed in our approach.

In the numerical simulation, we set $a_0 = 1.1196, k_{B}T_{0} = 0.25, k_{B}\Delta T = 0.05, m = 1.0, \gamma = 0.5$ and $N = 100$ in \eqref{eq:Lan_temp}. A Verlet-type of integration algorithm \cite{DBCC:98} is applied using step length $\Delta t = 0.05$. The data is sub-sampled at every ten integration time steps ($\delta t = 0.005$) at equilibrium with sample size $5\times 10^{6}$. The parameter $k_{B}T_{0}$ is set to be large enough, so that at the equilibrium, the marginal density of $\mbd$ is non-Gaussian. Throughout the section, the underlying dynamics \eqref{eq:Lan_temp} will only be used for generating CG observations and computing the full response for verification of the parameter estimates.

\subsection{The coarse-grained model}

Motivated by the residue-based coarse-graining 
\cite{baaden_coarse-grain_2013}, we define the CG variables by dividing the system equally to consecutive blocks, each of which contains $J$ atoms. In particular, the displacement and velocity of the CG particles are defined as
\begin{equation}\label{eq:CG_par}
q_{i} : = \frac{1}{J} \sum_{j = (i-1)J+1}^{iJ} r_j, \quad p_{i} := \dot{q}_{i} = \frac{1}{J} \sum_{j = (i-1)J+1}^{iJ} v_j, \quad i = 1,2, \dots, n_{b}, \quad n_{b}: = \frac{N}{J},
\end{equation}
respectively. Here, $n_{b}$ denotes the number of CG particles. To construct a CG model, we propose the following function of the CG relative displacement $\bm{\xi} = (q_{2}-q_{1}, q_{3}-q_{2}, \dots, q_{n_{b}}-q_{n_{b}-1}) \in \mathbb{R}^{n_b -1 }$ and CG velocity $\mbp$
\begin{eqnarray}\label{eq:CG_rho1}
   p^{\dagger}_{eq}(\bm{\xi}, \bm{p}) \approx \rho(\bm{\xi}, \bm{p}) \propto
    \prod_{i=1}^{n_b} \exp\left( - b_i \frac{p_{i}^2}{2} -\sum_{j=1}^4  a_{j,i}\xi_{i}^j -   \sum_{k=1}^3  c_{k,i}p_{i} \xi_i^{k} \right),\quad a_{4,i},\;b_{i}>0, \quad i = 1,2,\dots, n_b,
\end{eqnarray}
as the ansatz of the equilibrium distribution of the CG variables $(\bm{\xi}, \bm{p})$. In the exponent of Eq. \eqref{eq:CG_rho1}, the cross terms between $p_i$ and $\xi_i$ are introduced to incorporate the heat flux of the underlying dynamics \eqref{eq:Lan_temp} caused by the non-uniform temperature profile. When the temperature is uniform, there is no heat flux, and these cross terms would be zero (see Appendix~\ref{app:ansatz} for details).

Formally, to capture the proposed equilibrium distribution \eqref{eq:CG_rho1}, we define a pseudo-Hamiltonian
\begin{equation*}
    H:= \sum_{i=1}^{n_b} \left( b_{i}\frac{p_i^2}{2} + \sum_{j=1}^4 a_{j,i} \xi_{i}^{j} + \sum_{k=1}^{3} c_{k,i} p_{i}\xi^{k}_{i} \right),
\end{equation*}
such that the equilibrium density ansatz in \eqref{eq:CG_rho1} can be written as $\rho \propto \exp\left(-H\right)$. We consider the following imperfect model,
\begin{equation}\label{eq:CG_sgd}
\begin{split}
\dot{\mbq} &= D^{-1}\nabla_{\mbp}H, \\
\dot{\mbp} &= -D^{-1}\nabla_{\mbq}H - D^{-1} \Gamma \nabla_{\mbp}H + \sqrt{2} D^{-\frac{1}{2}} \sigma \dot{U}_t,
\end{split}
\end{equation}
where $D = \diag(b_1, b_2, \dots, b_{n_b})$, and $U_{t}$, to distinguish from $W_{t}$ in \eqref{eq:Lan_temp}, is an $n_{b}$-dimensional Wiener process. In \eqref{eq:CG_sgd}, the coefficient $\Gamma$ is assumed to be a symmetric positive definite (SPD) matrix satisfying $\Gamma = \sigma \sigma^{\top}$. In practice, $\sigma$ can be determined by a Cholesky decomposition of $\Gamma$. In particular, $\rho$ in \eqref{eq:CG_rho1} solves the stationary Fokker-Planck equation of \eqref{eq:CG_sgd}. This can be verified by direct computations.

To have a better understanding of the imperfect model \eqref{eq:CG_sgd}, we rewrite \eqref{eq:CG_sgd} as
\begin{equation}\label{eq:CG_sgd_per}
\begin{split}
\dot{\mbq} &= \mbp + f(\mbq), \\
\dot{\mbp} &= F(\mbq) - \tilde{\Gamma}\mbp + g
(\mbq, \mbp)  + \Sigma \dot{W}_t,
\end{split}
\end{equation}
where
\begin{equation*}
\begin{split}
\tilde{\Gamma} & = D^{-1}\Gamma D, \quad \Sigma = \sqrt{2}D^{-\frac{1}{2}}\sigma, \\
f_{i} & = \sum_{k=1}^{3} \frac{c_{k,i}}{b_{i}}\xi_{i}^{k}, \quad F_{i} = b_{i}^{-1}\sum_{j=1}^4 j \left(a_{j,i}\xi_{i}^{j-1} - a_{j,i-1}\xi_{i-1}^{j-1} \right), \\
g_{i} & = b_{i}^{-1} \sum_{k=1}^3 k \left(c_{k,i}p_{i}\xi_{i}^{k-1} - c_{k,i-1}p_{i-1}\xi_{i-1}^{k-1} \right), \quad i = 1,2,\dots, n_b.
\end{split}
\end{equation*}
Eq. \eqref{eq:CG_sgd_per} suggests that the imperfect model \eqref{eq:CG_sgd} is a perturbed dynamics of 
\begin{equation}\label{eq:Lan_non}
\begin{split}
\dot{\mbq} &= \mbp, \\
\dot{\mbp} &= F(\mbq) - \tilde{\Gamma}\mbp + \Sigma \dot{W}_t.
\end{split}
\end{equation}
Unlike the classical Langevin dynamics, the force $F$ in \eqref{eq:Lan_non} is nonconservative. In particular, 
\begin{equation*}
\frac{\partial F_{i}}{\partial q_{i+1}}  = b_{i}^{-1} \sum_{j=2}^{4}j(j-1)a_{j,i}\xi_{i}^{j-2} \not = 
\frac{\partial F_{i+1}}{\partial q_{i}}  = b_{i+1}^{-1} \sum_{j=2}^{4}j(j-1)a_{j,i}\xi_{i}^{j-2}.
\end{equation*}
By Theorem 3 in \cite{sachs2017langevin}, we can verify that \eqref{eq:Lan_non} is geometrically ergodic. Thus, such an imperfect model \eqref{eq:CG_sgd} can be interpreted as a perturbed dynamics of an ergodic Langevin dynamics with nonconservative forces. The ergodicity of the imperfect model \eqref{eq:CG_sgd} is left here as an open problem.

Motivated by the symmetry of the system, the damping matrix $\Gamma$ is assumed to be a symmetric circulant matrix:
\begin{equation}\label{eq:damping}
\Gamma = \begin{pmatrix}
\gamma_0 & \gamma_1 & \gamma_2  &  \cdots & \gamma_2 & \gamma_1 \\
\gamma_1 & \gamma_0 & \gamma_1 & \gamma_2 & \cdots & \gamma_2 \\
& \ddots & \ddots & \ddots &  &  \\
&   & \ddots & \ddots & \ddots &  \\
\gamma_2 & \cdots  & \gamma_2        & \gamma_1   & \gamma_0 &  \gamma_1 \\
\gamma_1 & \gamma_2 & \cdots & \gamma_2 & \gamma_1 & \gamma_0
\end{pmatrix} \in \mathbb{R}^{n_{b} \times n_{b}},
\end{equation}
which is fully determined by the first $\lfloor \frac{n_{b}}{2} \rfloor+1$ elements $(\gamma_{0},\cdots,\gamma_{\lfloor \frac{n_{b}}{2} \rfloor})$ in the first row. To guarantee that $\Gamma$ is positive definite, we pose a diagonal dominant constraint:
\begin{equation*}
\gamma_0 > 
\begin{cases}
& 2 \sum\limits_{i=1}^{\frac{n_{b}}{2}-1} \left|\gamma_i\right| + \left|\gamma_{\frac{n_{b}}{2}}\right|, \quad n_{b}\text{ even}, \\
& 2 \sum\limits_{i=1}^{\frac{n_{b}-1}{2}} \left|\gamma_i\right|, \quad n_{b}\text{ odd}.
\end{cases}
\end{equation*}
For simplicity, we will set $\gamma_{3}=\gamma_{4}=\cdots=\gamma_{\lfloor \frac{n_{b}}{2} \rfloor}=0$ in \eqref{eq:damping}, and the parameter space of the imperfect model \eqref{eq:CG_sgd} can be summarized as 
\begin{equation}\label{eq:Lan_para}
\theta:= \left\{b_{i}, a_{1,i},\ldots,a_{4,i}, c_{1,i},\ldots, c_{3,i}, \gamma_{0}, \gamma_{1},\gamma_{2}\right\}_{i=1}^{n_b}, \quad a_{4,i},\;b_{i}>0, \; i = 1,2,\dots, n_b, \quad \gamma_{0} > 2\left(\left|\gamma_{1}\right|+\left|\gamma_{2}\right|\right).
\end{equation}
To be consistent with the notation in Section~\ref{sec:scheme}, the parameters in \eqref{eq:Lan_para} include $\theta_{eq}:= \{b_i, a_{j,i}, c_{k,i}\}$, which are the coefficients in the density function \eqref{eq:CG_rho1}, and $\gamma_{i}$ are the component of the damping matrix \eqref{eq:damping}.

There are many other approaches in the model reduction of Langevin dynamics, e.g., the partitioning technique \cite{sweet2008normal}, the Petrov-Galerkin projection \cite{MLL:19}, and the maximum likelihood estimator \cite{dequidt2015bayesian}. Here, as we have emphasized, we consider the imperfect model \eqref{eq:CG_sgd} as a mean to introduce the model error, and use the one-dimensional system to demonstrate the parameter estimation procedure, instead of focusing on the model reduction technique.

\subsection{Parameter estimation and numerical results}

To this end, we have introduced our underlying dynamics \eqref{eq:Lan_temp} and the corresponding imperfect model \eqref{eq:CG_sgd} with unknown parameters $\theta$ in \eqref{eq:Lan_para}. In the numerical experiment, the observation of the CG model \eqref{eq:CG_sgd} will be acquired from Eq. \eqref{eq:CG_par} based on the time series of $(\mbd, \mbv)$ generated by the underlying dynamics. For instance, the CG relative displacement $\bm{\xi}$ can be determined by:
\begin{equation}\label{eq:CG_obs}
\xi_{i} = q_{i+1} - q_{i} = \frac{1}{J} \sum_{j = (i-1)J+1}^{iJ} (r_{j+J}-r_{j}) = \frac{1}{J} \sum_{j =(i-1)J+1}^{iJ} \;\;\sum_{k = j}^{j+J-1}(r_{k+1}-r_{k}), \quad i = 1,2,\dots, n_{b}-1.
\end{equation}
Here, Eq. \eqref{eq:CG_obs}, together with the second equation of  \eqref{eq:CG_par} correspond to the function $P(\cdot)$ in our general setup \eqref{eq:full_dyn} with the variable of interest $X=(\bm{\xi}, \mbp)$. Thus, we are able to obtain the time series of $(\bm{\xi}, \mbp)$, containing model error, for the imperfect model \eqref{eq:CG_sgd}, which will be the only information used in estimating $\theta$.

\begin{figure}[!ht]
	\centering
	\includegraphics[width=.49\textwidth]{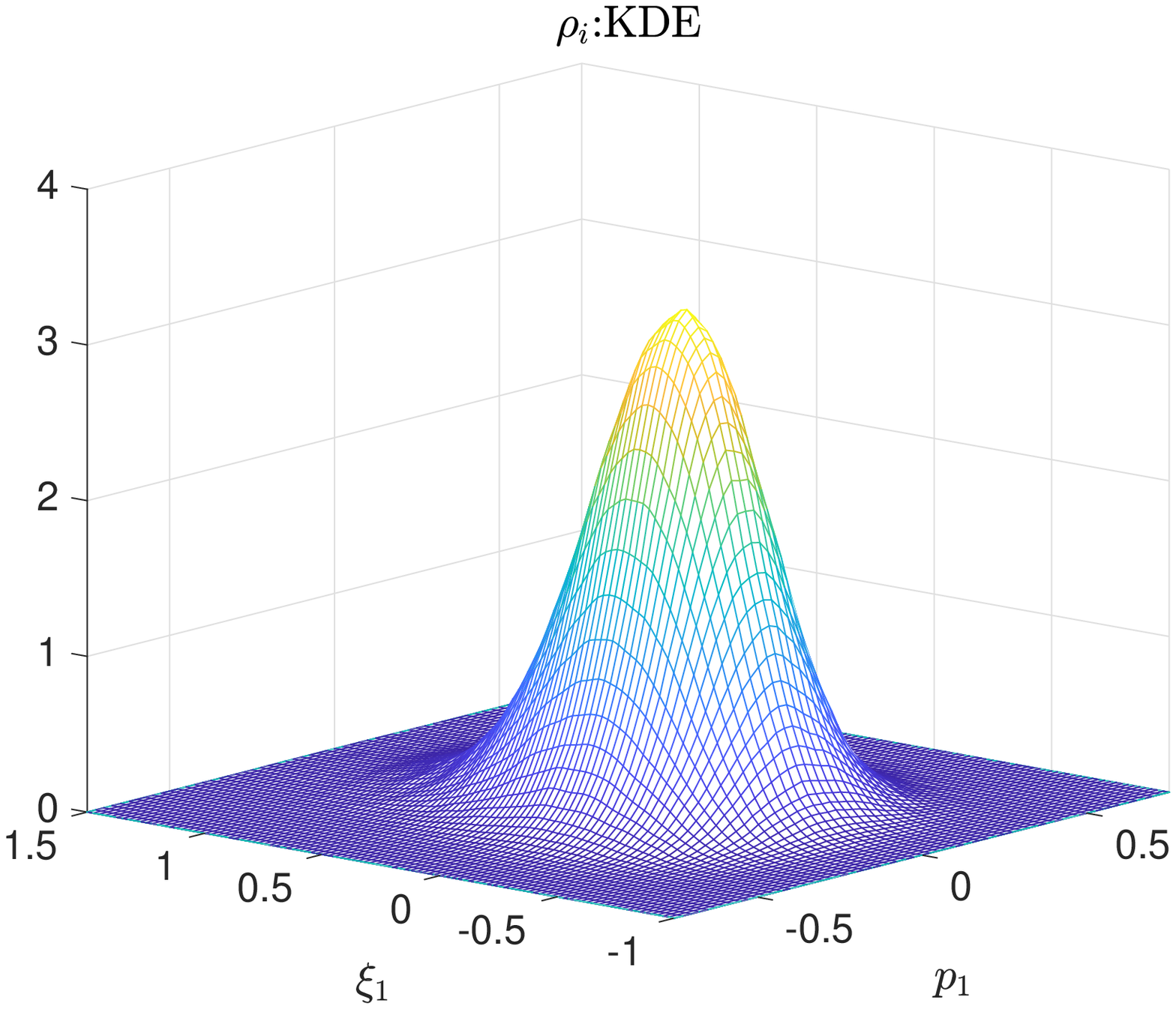}
	\includegraphics[width=.49\textwidth]{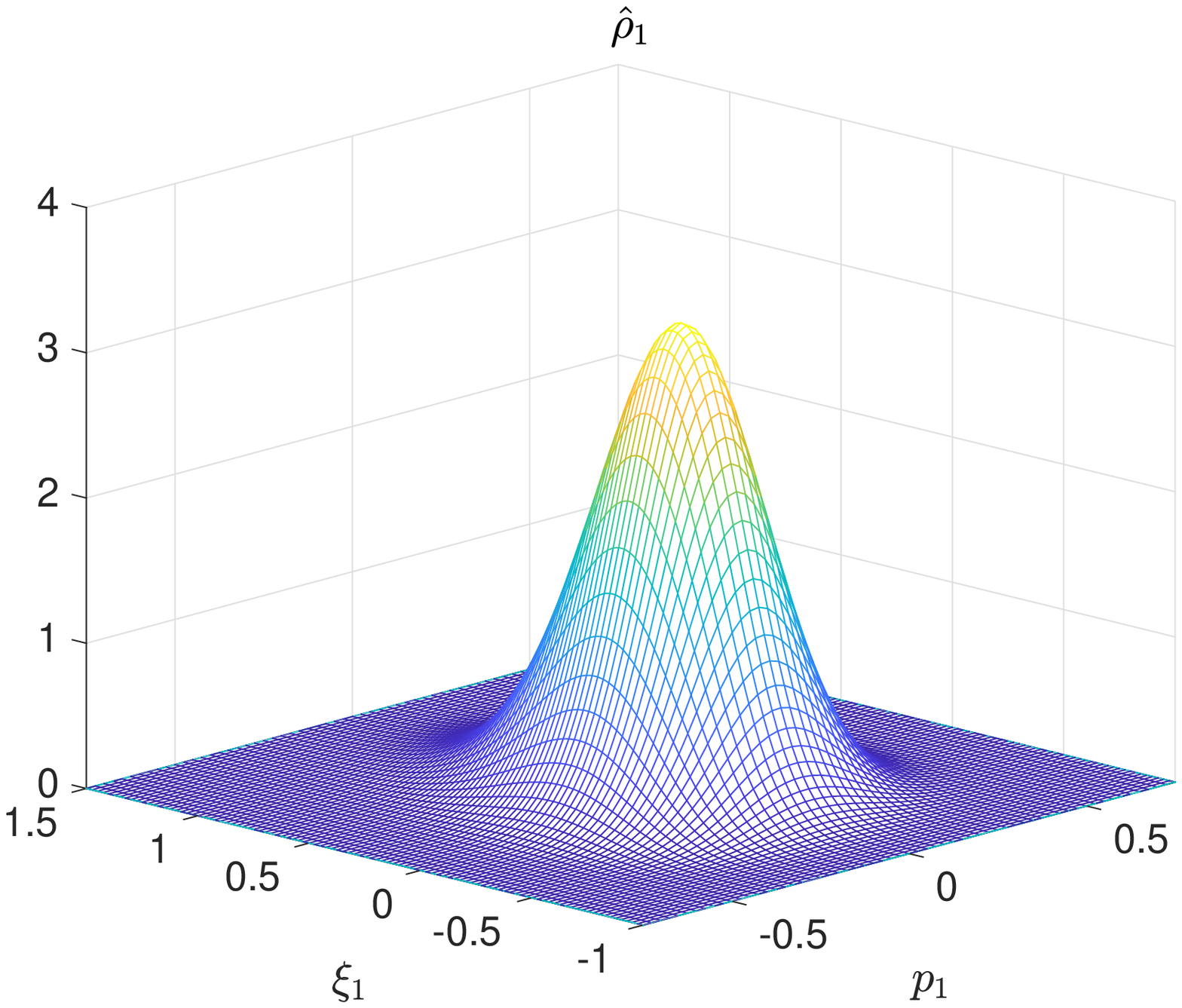}
	\includegraphics[width=.49\textwidth]{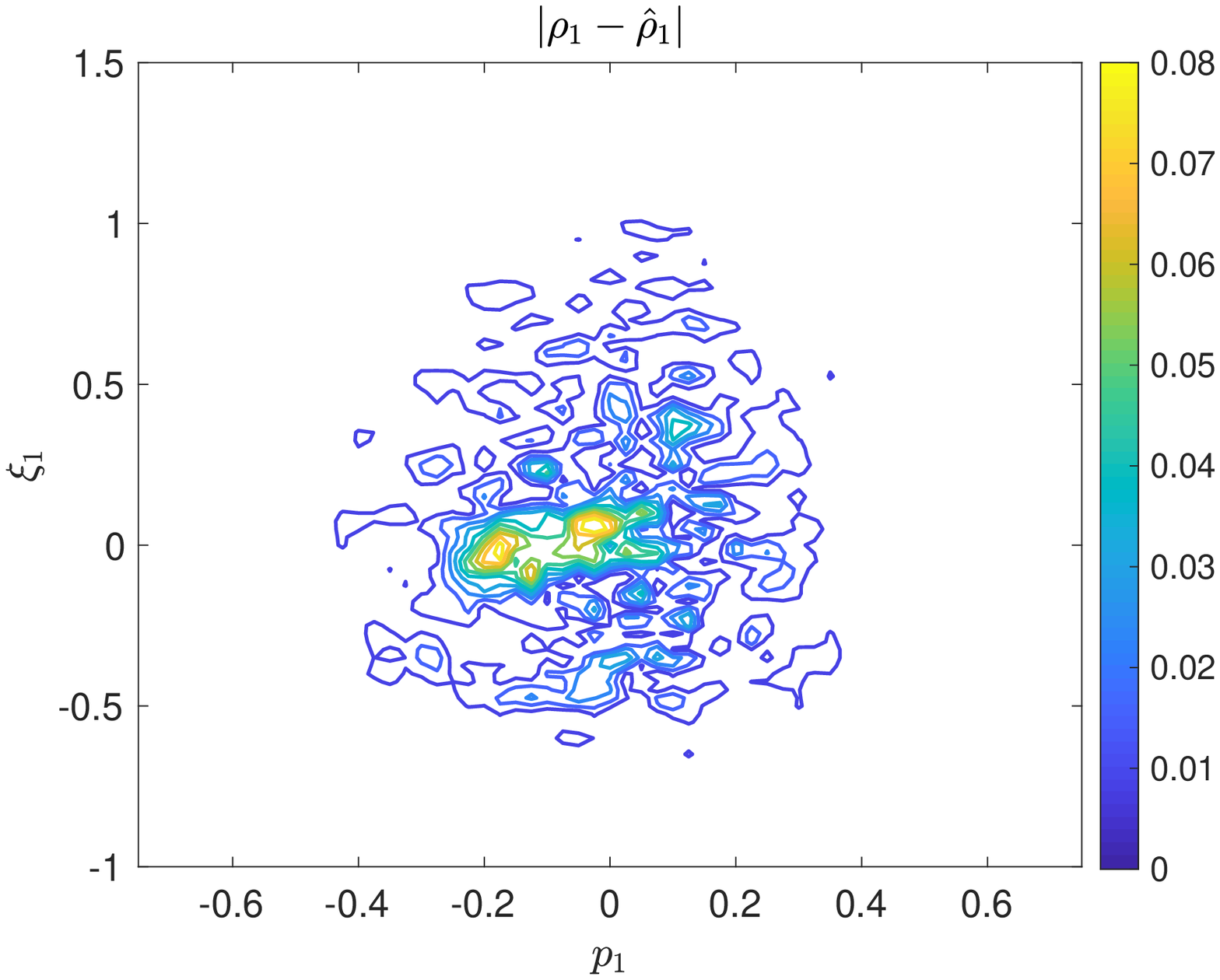}
	\caption{The density function $\rho_{1}$ (up-left, plot using the ``mvksdensity'' syntax in MATLAB), and $\hat{\rho}_1$ (up-right, plot using the parametric form in \eqref{eq:rho_i} and the parameter estimates in Table~\ref{tab:est_rho}) the density estimates generated via maximum entropy. The absolute error between $\rho_1$ and $\hat{\rho}_1$ is presented in the lower panel.}
	\label{fig:rho_1}
\end{figure}

For parameters $\theta_{eq}$ that appear in the density ansatz \eqref{eq:CG_rho1}, we will apply the maximum entropy method~\cite{jaynes:57} using the equation-by-equation solver introduced in \cite{hao2018equation}. In our implementation, we solve for $(b_i, a_{1,i},\ldots,a_{4,i}, c_{1,i},\ldots,c_{3,i})$ for each $i$ independently since
\begin{equation}\label{eq:rho_i}
   \rho_{i}(\xi_i, p_i) \propto \exp\left(- b_i \frac{p_{i}^2}{2} -\sum_{j=1}^4 a_{j,i}\xi_{i}^j - \sum_{k=1}^3  c_{k,i}p_{i} \xi_i^{k} \right)
\end{equation}
fully determines the marginal distributions of $(\xi_i,p_i)$. Thus, the parameter estimation problem of $\theta_{eq}$ can be decomposed into a total of $n_b$ two-dimensional maximum entropy problems. By maximizing the Shannon's entropy under the moment constraints, we formulate the following maximum entropy estimates for $i = 1,2,\dots, n_b$,
\BEA
&(\hat{b}_i, \hat{a}_{1,i},\ldots,\hat{a}_{4,i}, \hat{c}_{1,i},\ldots, \hat{c}_{3,i}) = \argmin \int -\rho_{i}\left(\xi_i, p_i \right) \log \rho_{i}\left(\xi_i, p_i \right) \td \xi_i \td p_i, \label{eq:est_a} \\
&\text{s.t. } \mathbb{E}_{\rho_i}[1] = 1,\; \mathbb{E}_{\rho_i}[p_i^2] = \mathbb{E}_{p^{\dagger}_{eq}} [p_i^2], \; \mathbb{E}_{\rho_i}[\xi_{i}^j] = \mathbb{E}_{p^{\dagger}_{eq}} [\xi_{i}^j], \; \mathbb{E}_{\rho_i}[p_{i} \xi_i^{k}] = \mathbb{E}_{p^{\dagger}_{eq}}[p_{i} \xi_i^{k}], \quad j = 1,2,3,4, \; k = 1,2,3, \label{eq:mom_con}
\EEA
where $\mathbb{E}_{p^{\dagger}_{eq}}[\cdot]$ will be approximated using Monte-Carlo average over the available time series of $(\bm{\xi}, \bm{p})$. Thus, such constrained minimization problem \eqref{eq:est_a} and \eqref{eq:mom_con} is equivalent to the reverse Monte Carlo method \cite{lyubartsev1995calculation}. Table~\ref{tab:est_rho} shows the value of the estimates of $\theta_{eq}$ for $i=1$. Figure~\ref{fig:rho_1} compares $\rho_{1}$ determined by a two-dimensional kernel density estimates based on the CG observations $(\xi_1, p_1)$ and $\hat{\rho}_1$ using the parameter estimates in Table~\ref{tab:est_rho}. We can see that the maximum entropy estimates provide a good fit towards the marginal distribution of $(\xi_1, p_1)$ at the equilibrium of the underlying dynamics.

\begin{table}[ht!]
    \centering
    \begin{tabular}{|c|c|c|c|c|c|c|c|c|}\hline
    Parameter ($\theta_{eq}$) & $b^{-1}_1$ & $a_{1,1}$ & $a_{2,1}$ & $a_{3,1}$ & $a_{4,1}$ & $c_{1,1}$ & $c_{2,1}$ & $c_{3,1}$ \\
    \hline
    Estimates ($\hat{\theta}_{eq}$)    & $ 0.02631 $ & $ -0.2582 $ & $ 5.951 $ & $ -2.781 $ & $ 0.6885 $ & $ -0.4949$ & $ 0.1446$ & $ 0.04487$ \\ \hline
    \end{tabular}
    \caption{Part of the estimates of $\theta_{eq}$ in \eqref{eq:CG_rho1}.}
    \label{tab:est_rho}
\end{table}

Up to this point, we have estimated $\theta_{eq}$ in the density ansatz \eqref{eq:CG_rho1} of the imperfect model \eqref{eq:CG_sgd}. We now apply our parameter estimation approach to estimate $\{\gamma_{i}\}$ in the damping matrix $\Gamma$. These parameters determine the average dynamical behavior at equilibrium. To satisfy the assumptions posed in Section~\ref{sec:scheme}, the external forcing of the underlying dynamics is set to be a constant external forcing, exerted equally on the first block of atoms, and the observables $A(X)$ are the velocities of the first two CG particles. As a result, following \eqref{eq:red_dyn_per}, the corresponding external forces $c\; \delta f$ and $c_{R}\;\delta f$ for the underlying dynamics (acting on variable $\bm{v}$) and imperfect model (acting on CG variable $\mbp$), respectively, satisfy
\begin{equation} \label{eq:extf_lan}
c = (\underbrace{1,\dots,1}_{J-\mbox{terms}}, 0, \dots, 0) \in \mathbb{R}^{N}, \quad 
c_{R} = (1,0,\dots, 0)^{\top}\in \mathbb{R}^{n_{b}}, \quad \delta f = \delta.
\end{equation}
According to the criterion postulated in Definition \ref{def:adm}, this forcing is admissible. With the observables $A = (p_{1},p_{2})^{\top}$, following \eqref{eq:lin_oper_red}-\eqref{eq:mar_esst}, we write
\begin{equation}\label{eq:Lan_esst}
\begin{aligned}
\hat{k}^{\dagger}_{A}(t) & =  \mathbb{E}_{p^{\dagger}_{eq}}\left[(p_{1}(t),p_{2}(t))^{\top}\otimes \hat{B}(p_{1}(0), \xi_1(0)) \right], \; \\
\hat k_{A}(t;\theta) &= \mathbb{E}_{\rho(\theta)}\left[(p_{1}(t),p_{2}(t))^{\top}\otimes \hat{B}(p_{1}(0), \xi_1(0)) \right],
\end{aligned}
\end{equation}
where 
\begin{equation*}
\hat{B}(p_{1}, \xi_1) = - \frac{\partial}{\partial p_{1}}\log\left( \rho(\bm{\xi}, \bm{p}; \hat{\theta}_{eq}) \right) =  \hat{b}_1 p_1 + \sum_{k=1}^{3}\hat{c}_{k,1} \xi_{1}^{k},
\end{equation*}
and $\hat{b}_{1}$ and $\{\hat{c}_{k,1}\}$ are the estimates obtained from solving the maximum entropy problem \eqref{eq:est_a}-\eqref{eq:mom_con}. Subsequently, we formulate the following constrained nonlinear least-squares problem:
\begin{equation} \label{eq:Lan_nls}
\min_{\theta} \sum_{i = 1}^{K} \sum_{j=1}^2 \left(g^{\dagger}_{j}(t_i) - \hat{g}_{j}(t_{i};\theta)\right)^2, \quad \text{s.t. } \theta_{eq} = \hat{\theta}_{eq}, \quad \gamma_{0} > 2\left(\left|\gamma_{1}\right|+\left|\gamma_{2}\right|\right),
\end{equation}
where
\begin{equation*}
g^{\dagger}_{j}(t) = \mathbb{E}_{p^{\dagger}_{eq}}\left[p_{j}(t)\hat{B}(p_{1}(0), \xi_1(0)) \right], \quad \hat{g}_{j}(t;\theta) = \mathbb{E}_{\rho(\theta)}\left[p_{j}(t)\hat{B}(p_{1}(0), \xi_1(0)) \right], \quad j = 1,2.
\end{equation*}
Recall that both $\{g^{\dagger}_{j}(t_i)\}$ and $\{\hat{g}_{j}(t;\theta)\}$ will be approximated by Monte-Carlo of the form \eqref{eq:MC_approx} based on, respectively, the available CG observations and the time series generated by simulating the imperfect model \eqref{eq:CG_sgd} at the corresponding parameter value. In the numerical simulation of the imperfect model \eqref{eq:CG_sgd}, we use a second-order weak trapezoidal method introduced in \cite{anderson2011weak}.

\begin{figure}[ht!]
	\begin{center}
		\includegraphics[width=.49\textwidth]{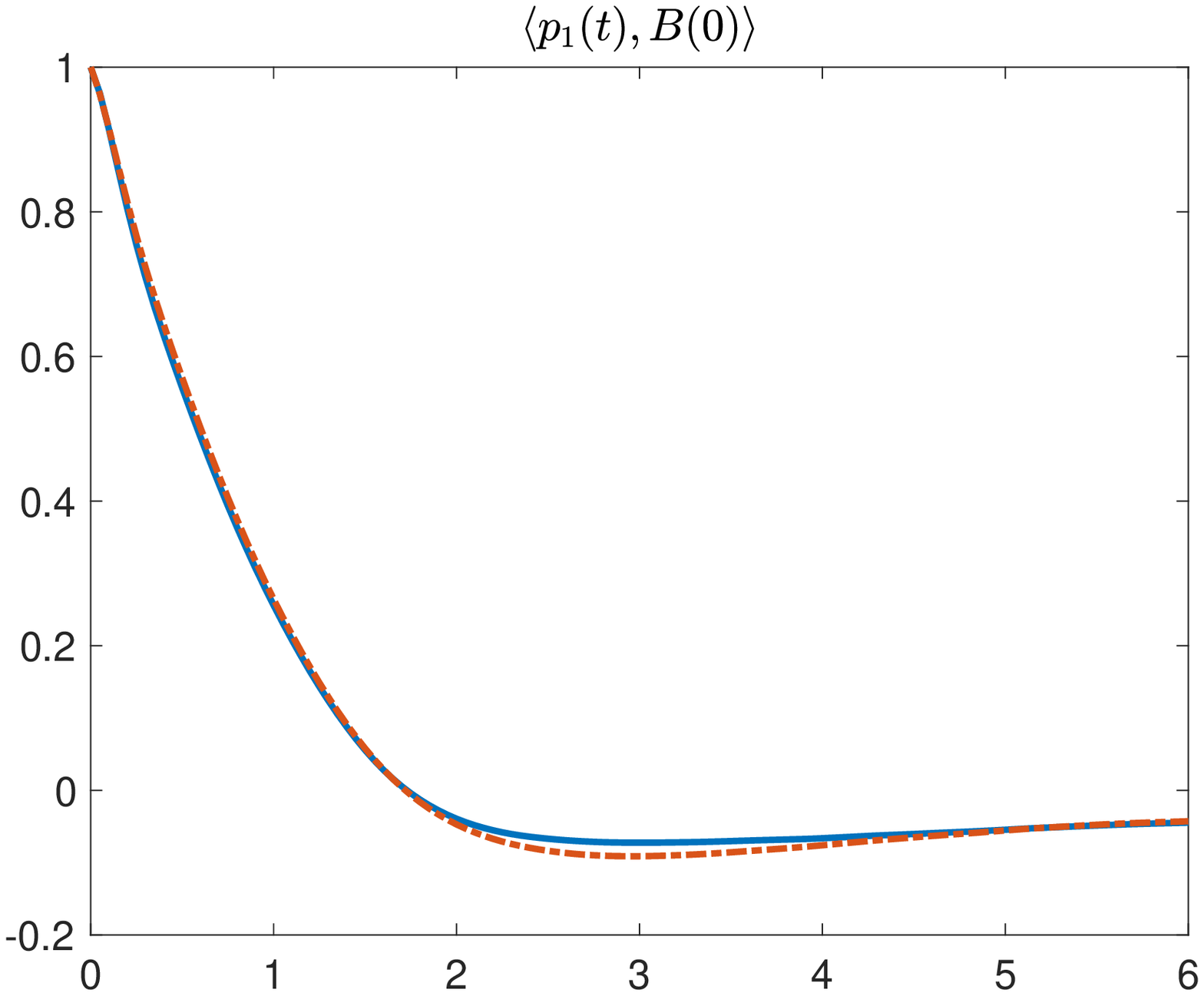}
		\includegraphics[width=.49\textwidth]{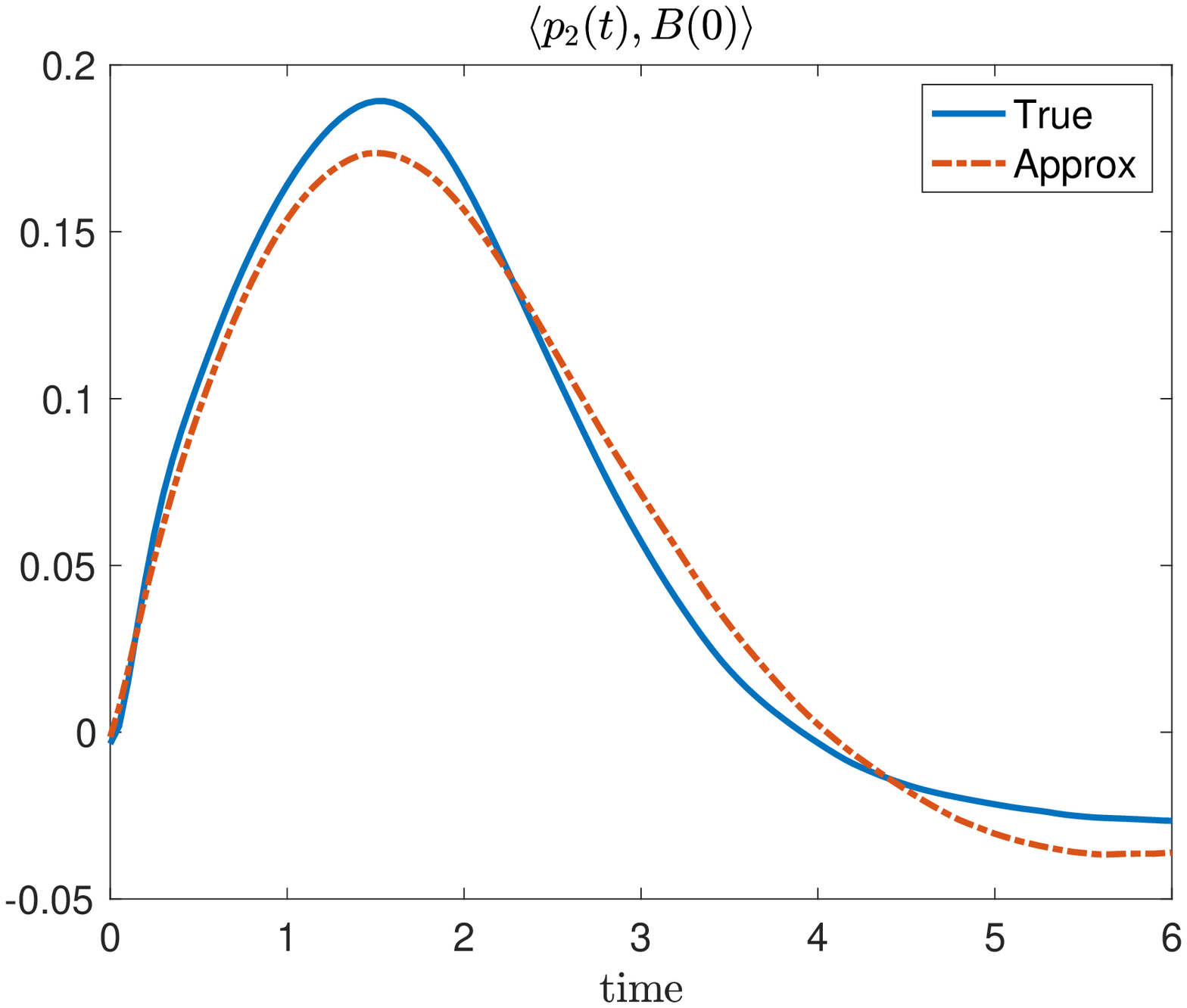}
		\caption{Comparison of the time correlation functions involved in the nonlinear least-squares problem \eqref{eq:Lan_nls}. The blue solids curves are computed from the available CG observations, while the red dot-dash curves are the approximation of the imperfect model based on the estimates. (sample size: $1\times 10^6$) }
		\label{fig:corr_fit}
	\end{center}
\end{figure}

\begin{figure}[ht!]
	\begin{center}
		\includegraphics[width=.49\textwidth]{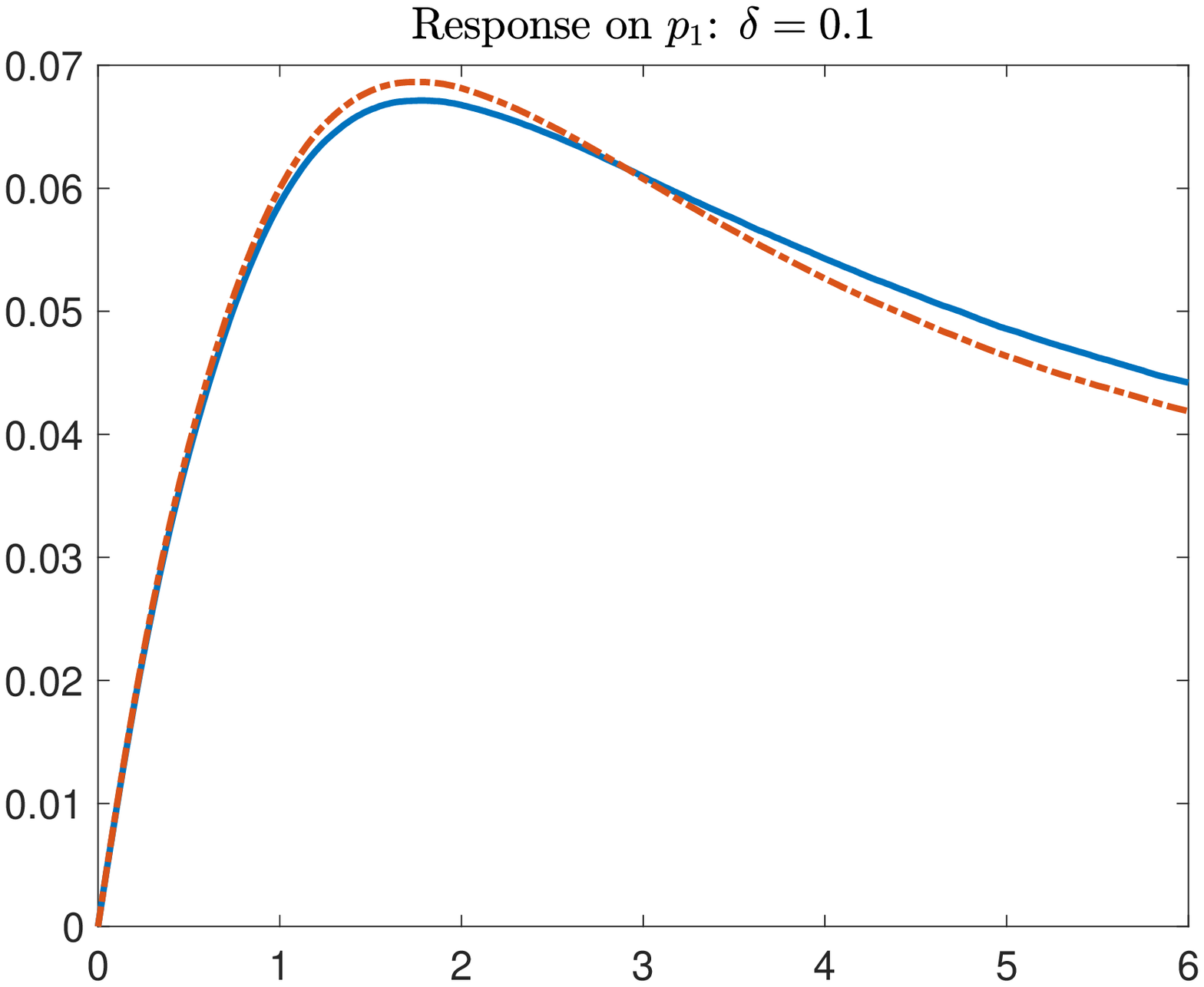}
		\includegraphics[width=.49\textwidth]{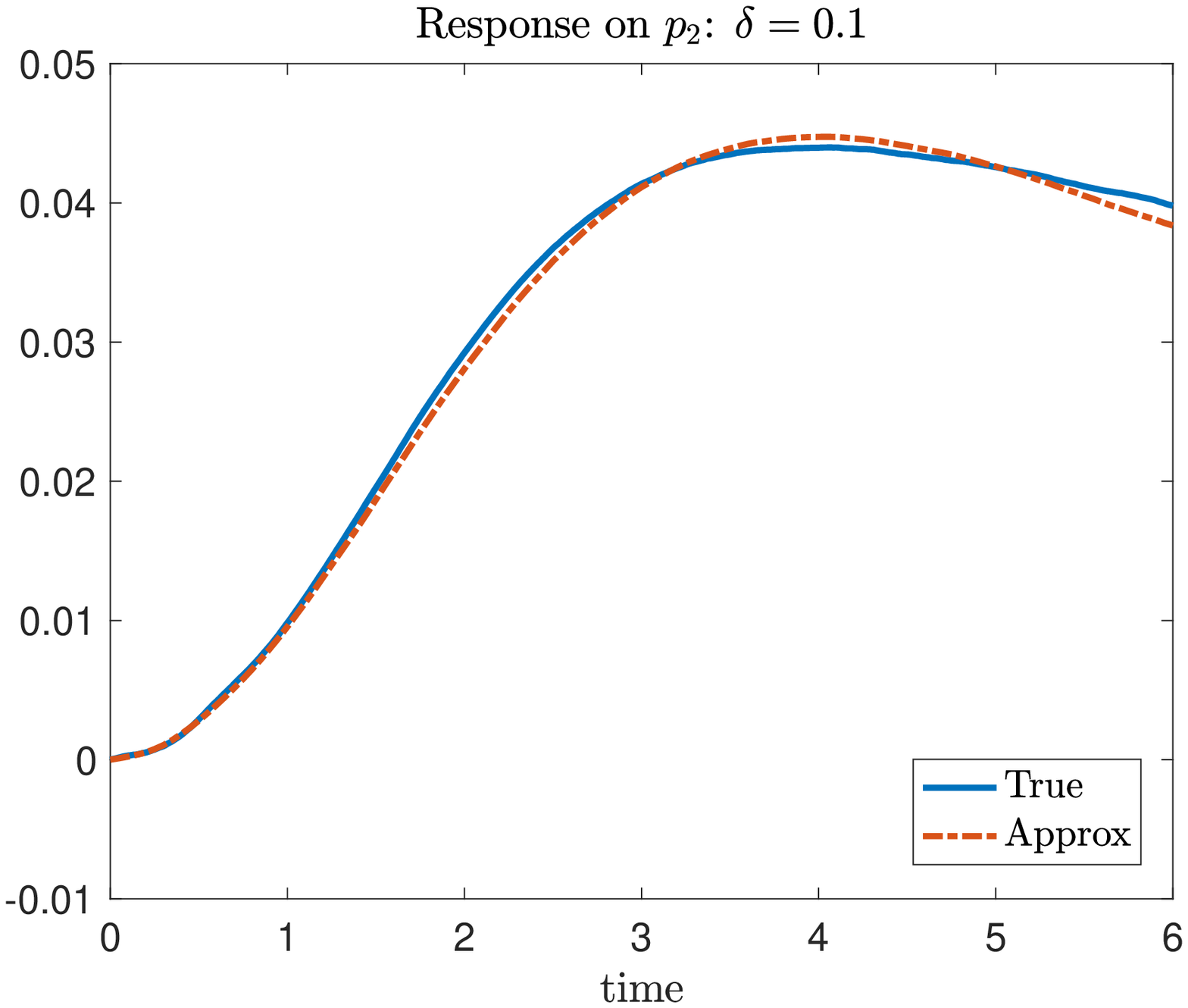}
		\caption{Reconstruction of full response statistics: underlying dynamics (blue solid) v.s. the imperfect model (red dash). We take $\delta = 0.1$ in \eqref{eq:extf_lan} as the external forcing. The full response statistics \eqref{eq:full_resp} are computed via Monte-Carlo. Check Appendix~\ref{app:A} for the details.}
		\label{fig:Lan_resp}
	\end{center}
\end{figure}

As an example, we set $J=10$ such that $n_{b}=10$ since $N=100$. Setting $t_i = 0.2i$ with $i = 1,2,\dots, 20$, we solve \eqref{eq:Lan_nls} using the numerical scheme reviewed in Section \ref{sec:review}. Specifically, a polynomial chaos expansion with order-$6$ Legendre polynomials is used to approximate the cost function \eqref{eq:Lan_nls} to avoid the excessive computational cost in evaluating $\hat{g}_j$ on new parameter $\theta$. (see \cite{HLZ:19} for details) Table \ref{tab:Lan_est} shows the value of the estimates. We also show in Figures~\ref{fig:corr_fit} and \ref{fig:Lan_resp} the results from fitting the essential statistics and recovering the full response of the underlying dynamics, respectively. One can observe from the figures that the estimates provide an excellent fit on the linear response operator $\langle p_1(t),B(0) \rangle$. For the response operator $\langle p_2(t),B(0) \rangle$, the poor performance at the longer time might be caused by the Monte-Carlo error involved in \eqref{eq:Lan_nls} and possibly the polynomial chaos approximation on the cost function. Interestingly, the imperfect model produces reasonably accurate estimates of the full nonlinear response statistics of the underlying dynamics. We also tested the isothermal case, as a special case of \eqref{eq:Lan_temp}, and the results also show good accuracy in terms of the estimation of the full response statistics. These  results are reported in Appendix~\ref{app:iso} for interested readers.

\begin{table}[ht]
	\begin{center}
		\begin{tabular}{|c|c|c|c|}
		\hline
			Parameter ($\gamma_i$) & $\gamma_{0}$ & $\gamma_{1}$ & $\gamma_{2}$ \\
			\hline
			Estimates ($\hat{\gamma}_i)$ & $0.8409$  & $-0.1758$ & $7.616\times 10^{-3}$ \\ \hline
		\end{tabular}
	\end{center}
    \caption{Estimates of the parameter in the damping matrix $\Gamma$ \eqref{eq:damping}.}
	\label{tab:Lan_est}
\end{table}

In the next section, we will present an example from solutions of a PDE. At this point, it is worthwhile to emphasize several points. 

\begin{rem}\label{rm:2}
While in general the choice of the imperfect model in \eqref{eq:red_dyn} depends on the specific applications, in current and next section, we will consider Langevin type stochastic dynamics, regardless of whether the underlying dynamics are Langevin or not. In the example above, the imperfect Langevin model \eqref{eq:CG_sgd}, as a result of coarse-graining, is a closed equation of $X = (\bm{\xi}, \mbp)$, with no dependence on any auxiliary variables $Y$. In this case, we propose a natural parametric form as an ansatz of the associated equilibrium density and use the maximum entropy principle to estimate $\theta_{eq}$. In the next section, we shall see that the imperfect model (also Langevin dynamics) is not a closed equation of $X$.
\end{rem}

\begin{rem}\label{rm:3} 
In the example above, we considered $p_R(X) = \rho(X;\theta_{eq})$ in \eqref{eq:CG_rho1}. We should point out that choosing $p_R$ to depend only on $X$ in \eqref{eq:mar_esst} is only sensible (in the linear response theory) if $X$ is independent to the orthogonal component, $X^\perp := u-P(u)$, at equilibrium state, in that $p_{eq}(X, X^{\perp}) = p_{R}(X)p_{R^{\perp}}(X^{\perp})$. In this case, we have
\begin{equation*}
\frac{\partial_{X_i}\left[c_{i}(X)p_{eq}(X, X^{\perp})\right]}{p_{eq}(X, X^{\perp})} = \frac{\partial_{X_i}\left[c_{i}(X)p_{R}(X)p_{R^{\perp}}(X^{\perp})\right]}{p_{R}(X)p_{R^{\perp}}(X^{\perp})}= \frac{\partial_{X_i}\left[c_{i}(X)p_{R}(X)\right]}{p_{R}(X)},
\end{equation*}
which suggests that the MLR operator is identical to the linear response operator of the underlying dynamics, $\hat{k}^{\dagger}_{A}(t) = k^{\dagger}_{A}(t)$. In general, however, the distributions of $X$ and $X^{\perp}$ are not independent and fitting to the marginal statistics defined through $p_R(X)$ may not be the optimal strategy. In the next section, we shall see that additional information is needed to achieve a reasonable linear response estimate.
\end{rem}

\section{A severely truncated spatio-temporal system} \label{sec:KS}

In this section, we consider a parameter estimation problem where the underlying equilibrium distribution is high-dimensional and unknown. With this constraint, there is a practical issue in the estimation of the linear response statistics $k^\dagger_A$ for solving the least-squares problem in \eqref{eq:nonlin_ls}, since the FDT approximation to these statistics requires an explicit form of the underlying density. To overcome this difficulty, we will consider the kernel embedding method \cite{Fuk:17,JH:18,JH:19,zhl:19b}  as a nonparametric estimator of the equilibrium density. We shall see that the resulting linear response statistics obtained from the kernel embedding estimate is more accurate than the Quasi-Gaussian FDT (QG-FDT) linear response \cite{mag:05,gritsun:07}. To test the performance of the parameter estimation method in recovering the full response statistics, we consider estimating parameters of an arbitrary closure model of a severely truncated Kuramoto-Shivashinsky (KS) equation, which full solutions exhibit spatiotemporal chaotic patterns observed in many applications, such as the trapped ion modes in plasma \cite{laquey1975nonlinear} and phase dynamics in reaction-diffusion systems \cite{10.1143/PTP.55.356}.  

\subsection{The underlying dynamics}

The underlying dynamics of the KS equation, defined on a one-dimensional periodic domain $[0,L]$, can be approximated on uniformly distributed spatial nodes (e.g.,\cite{Stinis:04,Trefethen:05,LLC:17}) with discrete Fourier modes that satisfy,
\begin{equation}\label{eq:Fourier_full}
\quad u_{0} = 0, \quad \dot{u}_{k} = (q_{k}^2 - q_{k}^{4})u_{k} - \frac{iq_{k}}{2} \sum_{1\leq |\ell|, |k-\ell|\leq \frac{N}{2}-1} u_{\ell} u_{k-\ell}, \quad |k|\leq \frac{N}{2}-1.
\end{equation}
Here, $q_{k} = \frac{2\pi k}{L}$ and $u_{k}(t)$ denote the Fourier coefficients that satisfy $u_{-k} = u_{k}^{*}$.

We take the numerical solution to the system \eqref{eq:Fourier_full} with $N=48$ as the solutions to the underlying dynamics; this choice of $N$ is sufficient to recover the inertial manifold of the system\cite{LLC:17}. The numerical approximate solution is generated using an exponential time difference fourth order Runge-Kutta (ETDRK4) method introduced in \cite{Trefethen:05}  with the discrete Fourier transform of $(1 + \sin(x))\cos(x)$ on the grid points $\{x_i=iL/N \}_{i=0}^{N-1}$  as the initial condition. Here, the ETDRK4 method is applied with $L = 2\pi / \sqrt{0.085}$. The value of $L$ is large enough to yield chaotic solutions \cite{Stinis:04}. We drop the first $5 \times 10^4$ time units to  allow the system to relax, and take $u_{1}$ and $u_{2}$ in the next $5 \times 10^5$ time units with step length $\Delta t = 0.1$ as the training data set.
We denote this training data set as $\{u_{1,t},u_{2,t}\}_{t=1}^{T}$ with $T=5 \times 10^6$.

We should point out that the chaotic behavior of the problem is characterized by positive Lyapunov exponents, the exponentially decaying time correlations, and other metrics \cite{toh:87} that are empirically estimated from the numerical solutions of \eqref{eq:Fourier_full}. Following many other works that used this model as a benchmark for recovery of the invariant statistics (e.g, \cite{Stinis:04,LLC:17,HJLY:19}), we assume that the numerical solutions of the KS system in \eqref{eq:Fourier_full} constitute a time series that is stationary in time and homogeneous in space. Effectively, this is the ergodicity assumption that is motivated by the empirical characterization of the chaotic behavior of the numerical solutions of \eqref{eq:Fourier_full} as discussed above.

\subsection{FDT response based on kernel embedding estimates} \label{sec:kme}

Our interest is in predicting the full response of the first Fourier mode, $u_1$, under small external forcing, without knowing the full underlying dynamics in \eqref{eq:Fourier_full}. Instead, we are only given the dynamical equation associated to the first Fourier mode,  
\begin{equation}\label{eq:missing}
\dot u_{1} = (q_{1}^2 - q_{1}^{4})u_{1} - iq_{1}u_{1}u_{2} + g, \quad g:=- \frac{iq_{k}}{2} \sum_{1\leq |\ell|, |k-\ell|\leq \frac{N}{2}-1} u_{\ell} u_{k-\ell}+iq_{1}u_{1}u_{2} 
\end{equation}
where $g$ is the {\it identifiable unresolved component}
in the following sense. Given the partial dynamics \eqref{eq:missing} and the training data set $\{u_{1,t},u_{2,t}\}_{t=1}^{T}$, one can extract $\{g_t\}_{t=1}^{T}$ by fitting these data set into the first equation in \eqref{eq:missing}, where the time derivative is replaced by a finite difference approximation \cite{LLC:17}.
 
\begin{figure}[ht]
	\begin{center}
		\includegraphics[width=.49\textwidth]{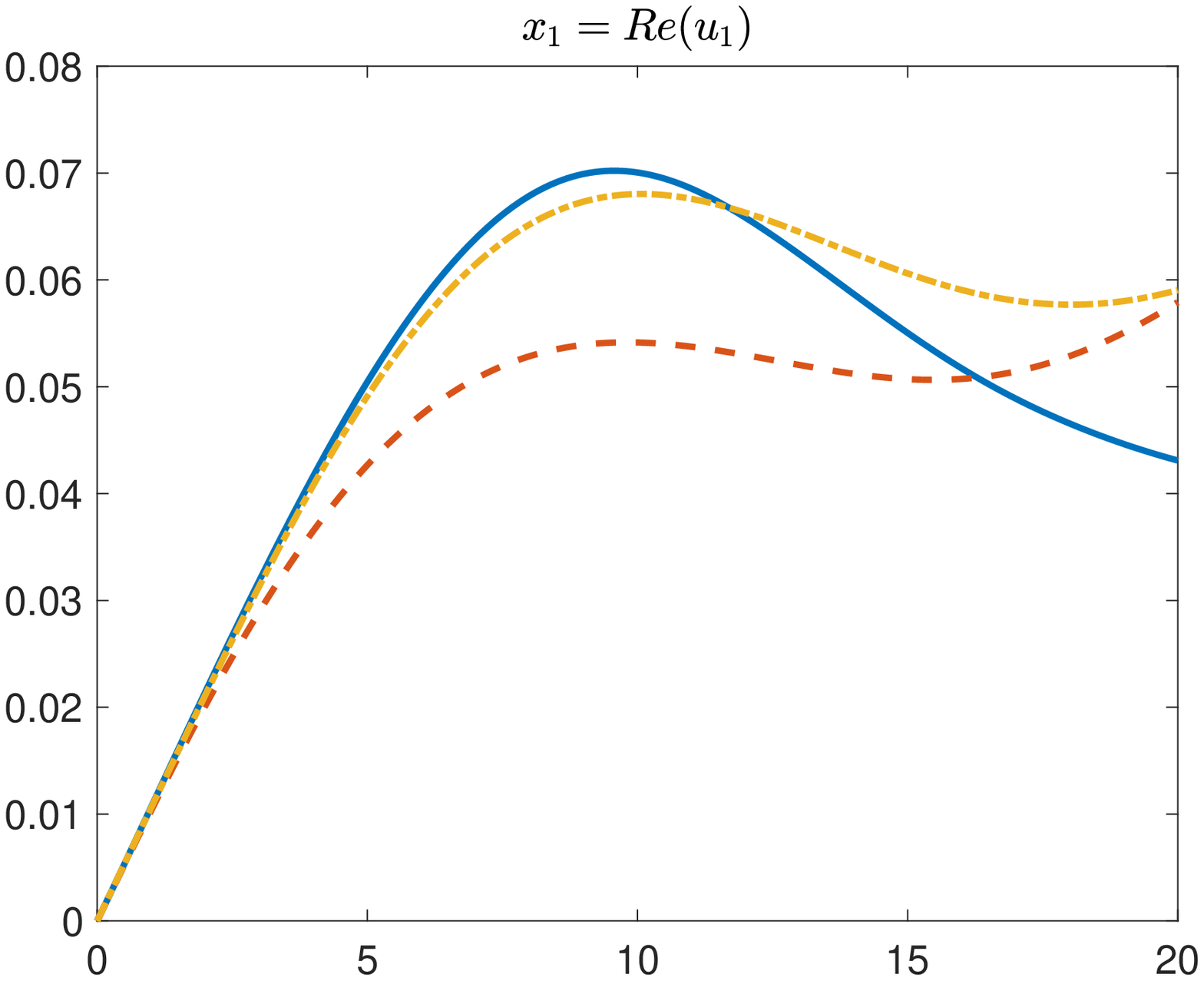}
		\includegraphics[width=.49\textwidth]{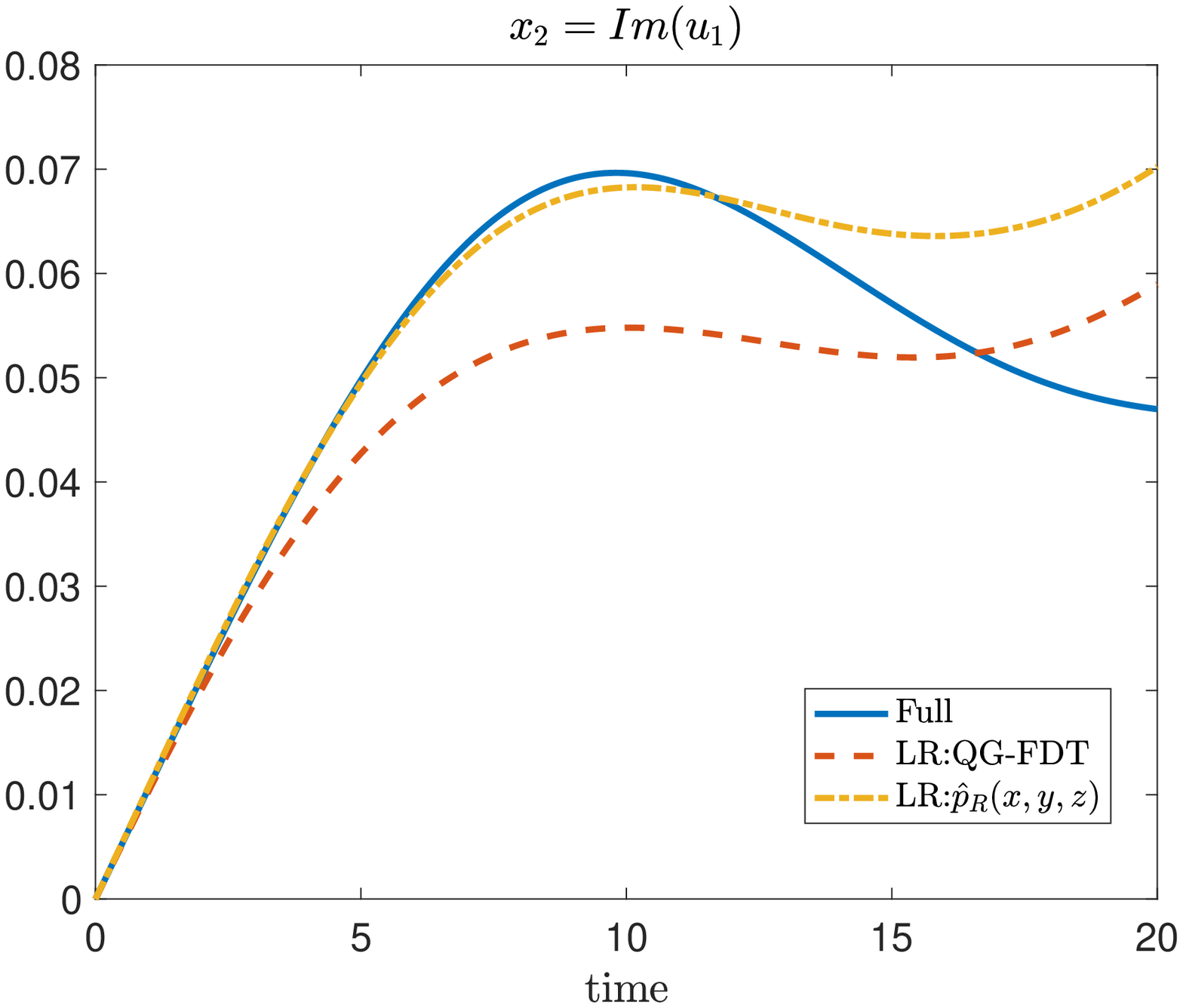}
		\caption{The response on $u_{1}$ component of the underlying dynamics \eqref{eq:Fourier_full} to a constant external forcing $\delta f = \frac{1}{N}(0.5 + 0.5 i )$. The blue solid curves show the full nonlinear response of $u_1$ approximated by Monte-Carlo based on $5\times 10^{5}$ times of simulations; the red dash shows the performance of QG-FDT based on the first six modes $(u_{1},u_{2}, \dots, u_{6})$; yellow dot-dash curves correspond to the MLR statistics produced by \eqref{eq:lin_resp_KS1} based on ($u_{1}, u_{2}, g$) as the resolved variable.}
		\label{fig:KS_1}
	\end{center}
\end{figure}

As we previously alluded, the key issue here is that we have no access to the explicit form of the equilibrium distribution of the full dynamics, which is very high-dimensional. Since we only have data set $\{u_{1,t},u_{2,t},g_t\}_{t=1}^{T}$, we propose to consider the linear response operator in \eqref{eq:lin_oper}, computed by averaging over $B$ that is defined with respect to a marginal density, $p_R$, of these three variables $(u_1,u_2,g)$. Beyond this practical consideration, we should allude that the inclusion of the identifiable unresolved variable $g$ is essential. To understand this, we first point out that if we follow the notation in \eqref{eq:full_dyn}, then $X$ corresponds to $u_1$. If we only consider the MLR operator with $p_R$ that depends only on $u_1$ (or even on $u_1,u_2$), we found that the corresponding linear response estimates do not approximate well the full response statistics. This is not surprising since $X = u_1$ is not independent to the conjugate variables, $X^\perp = (u_2,\ldots,u_{N})$, as discussed in Remark~\ref{rm:3}. Likewise, $(u_1,u_2)$ is not independent to its conjugate variables $(u_3,\ldots,u_{N})$. As we shall see below (see Figure~\ref{fig:KS_1}), the resulting MLR statistics corresponds to $p_R$ that depends on $(u_1,u_2,g)$, is an accurate estimator of the full nonlinear response statistic of the underlying dynamics under admissible external disturbances up to a finite time. The numerical results suggest that the inclusion of $(u_2,g)$ in $p_R$ is important since it carries information in $X^\perp$ that is relevant for accurate linear response estimations.

To facilitate the following discussion, we define
\begin{equation}\label{eq:real_rep}
x: = \big(\opre(u_{1}), \opim(u_{1})\big)^{\top}, \quad 
y: = \big(\opre(u_{2}), \opim(u_{2})\big)^{\top}, \quad \text{and} \quad
z: = \big(\opre(g), \opim(g)\big)^{\top}.
\end{equation}
Under a constant external forcing on $u_{1}$, which is admissible according to the criteria in Definition~\ref{def:adm}, Eq. \eqref{eq:mar_esst} suggests that, using the real representation \eqref{eq:real_rep}, the MLR  operator of \eqref{eq:Fourier_full}  is given by
\begin{equation}\label{eq:lin_resp_KS}
\hat{k}^{\dagger}_{A}(t) = \mathbb{E}_{p_{eq}^\dagger}\left [A(x(t))\otimes B(x(0), y(0), z(0)) \right], \quad B(x,y,z) := -\nabla_{x}\log\left(p_R(x,y,z)\right),
\end{equation}
where the unknown $p_R$ corresponds to the marginal distribution of $(x,y,z)$ at the equilibrium of the underlying dynamics \eqref{eq:Fourier_full}. While the expectation over $\peq^\dagger$ can be realized using Monte-Carlo average over the available data, one needs an actual expression of $p_R$ which appears in $B$. This way, we have replaced the problem of estimating a very high-dimensional density of $(u_1,\ldots,u_N)$ to a moderately low-dimensional (in this case, six-dimensional real-valued) function. In the remainder of this section, we will discuss how to approximate $p_R$ with a kernel embedding estimator, $\hat{p}_R$, and check the validity of the MLR in predicting the full response statistics of the underlying dynamics under small perturbations.

Recall that kernel embedding, as a nonparametric approach to learn a distribution function from the sample, can be formulated with an arbitrary kernel that uniquely determines the reproducing kernel Hilbert space (RKHS). We consider kernels defined by the orthogonal basis of a weighted $L^2$-space. Specifically, let the target density function $f \in \mathcal{H}_{\rho} \subset L^{2}\left(\mathbb{R}^{d}, q^{-1}\right)$, where $\mathcal{H}_{\rho}$ denotes the RKHS induced by a set of orthonormal basis $\left\{\Psi_{\vec{m}} = \psi_{\vec{m}} q\right\}$ of $L^{2}\left(\mathbb{R}^{d}, q^{-1}\right)$ and $q$ denotes a positive weight function on $\mathbb{R}^{d}$, (see \cite{JH:18,JH:19} for details) where $\vec{m} = (m_{1}, m_{2}, \dots, m_{d})$ is a multi-index notation with nonnegative components.  In particular, we take $q$ to be the standard $d$-dimensional Gaussian distribution, and define the kernel
\begin{equation} \label{eq:kernel}
k_{\rho}(\bm{x}, \bm{y}): = \sum_{\vec{m}\geq 0} \lambda_{\vec{m}} \Psi_{\vec{m}}(\bm{x}) \Psi_{\vec{m}}(\bm{y}), \quad \lambda_{\vec{m}}: = \rho^{\|\vec{m}\|_{1}}, \quad \|\vec{m}\|_{1}: = \sum_{k=1}^{d} m_{k},
\end{equation}
for $\rho\in (0,1)$. In \cite{zhl:19b}, we have shown that equipped with the inner product
\begin{equation} \label{eq:inner_prod}
 \langle f, g\rangle_{\mathcal{H}_{\rho}} : = \sum_{\vec{m}\geq 0} \frac{\hat{f}_{\vec{m}}\hat{g}_{\vec{m}}}{\rho^{\|\vec{m}\|_{1}}}, \quad f = \sum_{\vec{m}\geq 0} \hat{f}_{\vec{m}}\Psi_{\vec{m}}, \; g = \sum_{\vec{m}\geq 0} \hat{g}_{\vec{m}}\Psi_{\vec{m}},
\end{equation}
$\mathcal{H}_{\rho}$ is an RKHS with a bounded kernel $k_{\rho}$ satisfying
\begin{equation*}
k_{\rho}(\bm{x}, \bm{y})  = (2\pi)^{-d}\left(1-\rho^2\right)^{-\frac{d}{2}} \exp \left[ -\frac{1}{2(1-\rho^2)}\left(\|\bm x\|^2 + \|\bm y\|^2 - 2\rho \sum_{i=1}^{d} x_{i}y_{i} \right)   \right].
\end{equation*}
Here, $k_{\rho}$ \eqref{eq:kernel} can be interpreted as a generalization of Mehler kernel \cite{mehler1866ueber}. Then the order-$M$ kernel embedding estimate of $f$, denoted by $f_{M}$, is given by
\begin{equation}\label{eq:kme1}
f_{M} := \sum_{\|\vec{m}\|_{1}\leq M} \hat{f}_{\vec{m}} \Psi_{\vec{m}}, \quad \hat{f}_{\vec{m}} = \int_{\mathbb{R}^{d}} f\Psi_{\vec{m}} q^{-1} \td \bm{x} = \int_{\mathbb{R}^{d}} f\psi_{\vec{m}}\td \bm{x}.
\end{equation}
Since $\{\psi_{\vec{m}}\}$ forms an orthonormal basis for $L^{2}(q)$, we take
\begin{equation*}
\psi_{\vec{m}}(\bm{x}) = \prod_{k = 1}^{d} \psi_{m_{k}}(x_{k}), \quad \bm{x} = (x_{1},x_{2}, \dots, x_{d}),
\end{equation*}
where $\psi_{m_{k}}$ is the order-$m_{k}$ normalized Hermite polynomial. We should point out with this choice of basis representation, we basically arrive at a polynomial chaos approximation \cite{xiu:2010} of a density function $f$. In practice, the integral in \eqref{eq:kme1} can be approximated by Monte-Carlo,
\begin{equation}\label{eq:kme_MC}
\hat{f}_{\vec{m}} \approx \hat{f}_{\vec{m},N} =  \frac{1}{N} \sum_{n=1}^{N} \psi_{\vec{m}}(X_{n}),
\end{equation}
where $\{X_{n}\}_{n=1}^{N}$ are sampled from the target density function $f$. One can define the order-$M$ empirical kernel embedding estimate of $f$ as
\begin{equation}\label{eq:kme2}
f_{M,N} : = \sum_{\|\vec{m}\|_{1}\leq M} \hat{f}_{\vec{m},N} \Psi_{\vec{m}}.\end{equation}
In our application, the sample $\{X_{n}\}$, as a time series of $(x,y,z)$ \eqref{eq:real_rep}, can be interpreted as a stationary process with the target distribution function $f = p_R$. 

Take the observable $A(x,y,z)=x$, and let $\hat{p}_{R,M}$ be the order-$M$ kernel embedding estimator of $p_R$. Then, Eq. \eqref{eq:lin_resp_KS} is approximated by an order-$M$ MLR estimator,
\begin{equation}\label{eq:lin_resp_KS1}
\hat{k}^{\dagger}_{A}(t) \approx -\mathbb{E}_{p_{eq}^\dagger}\left[x(t)\otimes \nabla_{x}\log(\hat{p}_{R,M}(x(0),y(0),z(0)))\right].
\end{equation}
Recall that the expectation is still over $p_{eq}^\dagger$ since we are using the observed time series from the underlying dynamics to compute the statistics \eqref{eq:lin_resp_KS1}. We should point out that $\nabla_x \log(\hat{p}_{R,M})$ in Eq. \eqref{eq:lin_resp_KS1} is not an accurate \emph{pointwise} estimator for $\nabla_x \log(p_{R})$ since approximating $\nabla_x \log(p_{R})$ is numerically ill-posed when $p_{R}$ is close to zero. However, estimating \eqref{eq:lin_resp_KS} with the estimator \eqref{eq:lin_resp_KS1} is not an ill-posed problem since the MLR operator $\hat{k}^{\dagger}_{A}$, as a statistical quantity, is defined as an expectation with respect to $p_R$, which cancels out the pointwise error.

\begin{rem}
Theoretically, the two-point statistics on the right-hand side of \eqref{eq:lin_resp_KS1} may not be well-defined, since the basis function $\Psi_{\vec{m}}$ used in the order-M kernel embedding estimate \eqref{eq:kme1} can be negative. We resolve this issue by restricting $\hat{p}_{R,M}$ in \eqref{eq:lin_resp_KS1} to $D_{M,\delta} = \{\hat{p}_{R,M}\geq \delta \}$ for $\delta>0$ such that the two-point statistics in \eqref{eq:lin_resp_KS1} is well-defined and the error is negligible \cite{zhl:19b}. In practice, fix a $\delta>0$, for large enough $M$, the sample points almost surely stay inside $D_{M,\delta}$, and the MLR estimator can be approximated by a Monte-Carlo integral over the entire data set.
\end{rem}

Using the approximate response operator in \eqref{eq:lin_resp_KS1}, we now compare the resulting approximate linear response statistics with the corresponding full response statistics. The numerical result (Figure~\ref{fig:KS_1}) suggests that the full response is well captured by the linear response that uses $\hat{p}_{R,M}(x,y,z)$ up to finite time. This result empirically validates the MLR response in this model. In fact, the MLR statistic is more accurate compared to the QG-FDT response \cite{mag:05,gritsun:07}, which is obtained from the response operator \eqref{eq:MC_approx} with $B^\dagger$ replaced by the one defined with a Gaussian approximation of the invariant density.

\subsection{The truncated model}\label{sec:reduce_KS}

The focus in this section is to devise a simple dynamical model that can improve the estimation of the full response, yet retaining the accuracy of the equilibrium distribution of variable $u_1$ obtained through the kernel embedding approach discussed in Section~\ref{sec:kme}. We should point out that while many reduced-order models \cite{Stinis:04,Pavliotis:13,LLC:17,HJLY:19} (just to name a few) have been proposed for the KS model, their constructions typically require the knowledge of more than two Fourier modes. In our case, we only work with the data set of the first two Fourier modes and the resulting identifiable irrelevant variables, $g$, as defined in \eqref{eq:missing}. While the dynamical equation for the first component in \eqref{eq:missing} is given, the fact that the linear term is unstable, since $q_1^2-q_1^4 >0$, makes it difficult to devise a closure model that is stable and ergodic with accurate invariant density of the variable $x: = \left(\opre(u_{1}), \opim(u_{1})\right)^{\top}$.

Given such constraints, we propose the following ``semi-parametric'' extended Langevin equation as our imperfect model
\begin{equation}\label{eq:KS_red}
\begin{split}
\dot x &= v, \\
\dot v &= \Lambda \nabla_x \log\left( \hat{\rho}(x)\right) - \Gamma v + \sigma \dot{W}_{t},
\end{split}
\end{equation}
where $x\in \mathbb{R}^2$ is the real representation of $u_{1}$ as in \eqref{eq:real_rep} and $W_{t}$ denotes the standard two-dimensional Wiener process. By ``semi-parametric'', we refer to the combination of linear parametric equation in the right-hand-side with a ``nonparametric'' term that involves  $\hat{\rho}(x)$ that is estimated by the kernel embedding approximation of the marginal distribution of $x$ at equilibrium. We should point out that since the density $\hat{\rho}$ is approximated using a specific choice of kernel with Hermite polynomials, as explained in Section~\ref{sec:kme}, the resulting model in \eqref{eq:KS_red} is effectively parametric.

\begin{rem}
Assume the extended Langevin equation \eqref{eq:KS_red} satisfies the following conditions.
\begin{enumerate}
    \item $\Lambda$ is SPD, $\Gamma$ is positive definite, and $\sigma$ satisfies the Lyapunov equation
    \begin{equation} \label{eq:lya_eq}
    \Gamma\Lambda + \Lambda \Gamma^{\top} = \sigma\sigma^{\top}.
    \end{equation}
    \item $\hat{\rho}(x)$ is a bounded smooth density function defined on $\mathbb{R}^{2}$, and there exist constants $D,E>0$ and $F$, such that
    \begin{equation*}
 -\langle \nabla_{x} \log (\hat{\rho}(x)) , x \rangle \geq -D \log(\hat{\rho}(x)) + E \|x\|_{2}^{2} + F, \quad \forall x \in \mathbb{R}^2.
\end{equation*}
\end{enumerate}
Then the geometric ergodicity of the imperfect model \eqref{eq:KS_red} follows from Theorem 3.2 in \cite{Mattingly:02}. In particular, consider a change of variables
\begin{equation*}
q := \Lambda^{-\frac{1}{2}}x, \quad p := \Lambda^{-\frac{1}{2}}v,
\end{equation*}
then $(q, p)$ satisfies 
\begin{equation} \label{eq:Std_Lan}
\begin{split}
\dot{q} &= p, \\
\dot{p} &= \nabla_{q} \log \left(\hat{\rho}(\Lambda^{\frac{1}{2}}q) \right) - \tilde{\Gamma}p + \tilde{\sigma} \dot{W}_{t},
\end{split}
\end{equation}
with $\tilde{\Gamma} = \Lambda^{-\frac{1}{2}} \Gamma \Lambda^{\frac{1}{2}}$, $\tilde \sigma = \Lambda^{-\frac{1}{2}} \sigma$, which suggests that 
\begin{equation}\label{eq:peq_red_KS}
    \rho_{eq}(x,v) \propto \hat{\rho}(x) \exp\left(-\frac{1}{2}v^{\top} \Gamma^{-1} v \right)
\end{equation}
is the unique equilibrium distribution of the imperfect model \eqref{eq:KS_red}.

In this example, choosing Langevin dynamic as the imperfect model is motivated largely by the separable structure of its equilibrium distribution as shown in \eqref{eq:peq_red_KS}. Following the notation in \eqref{eq:lin_oper_red}, $X=x$ and $Y=v$. In this case, the auxiliary variable $v$ is Gaussian with parameters that can be estimated directly from its empirical variance. As we have mentioned at the beginning of this section, while other choices of imperfect models that involve variables $(u_2, g)$ with accurate equilibrium density estimate $\hat{\rho}(x)$ are possible, they tend to be more complicated and required higher wavenumbers ($u_3, u_4, u_5, u_6$) as reported in \cite{LLC:17,HJLY:19}.
\end{rem}

\begin{figure}[ht]
	\begin{center}
	\includegraphics[width=.49\textwidth]{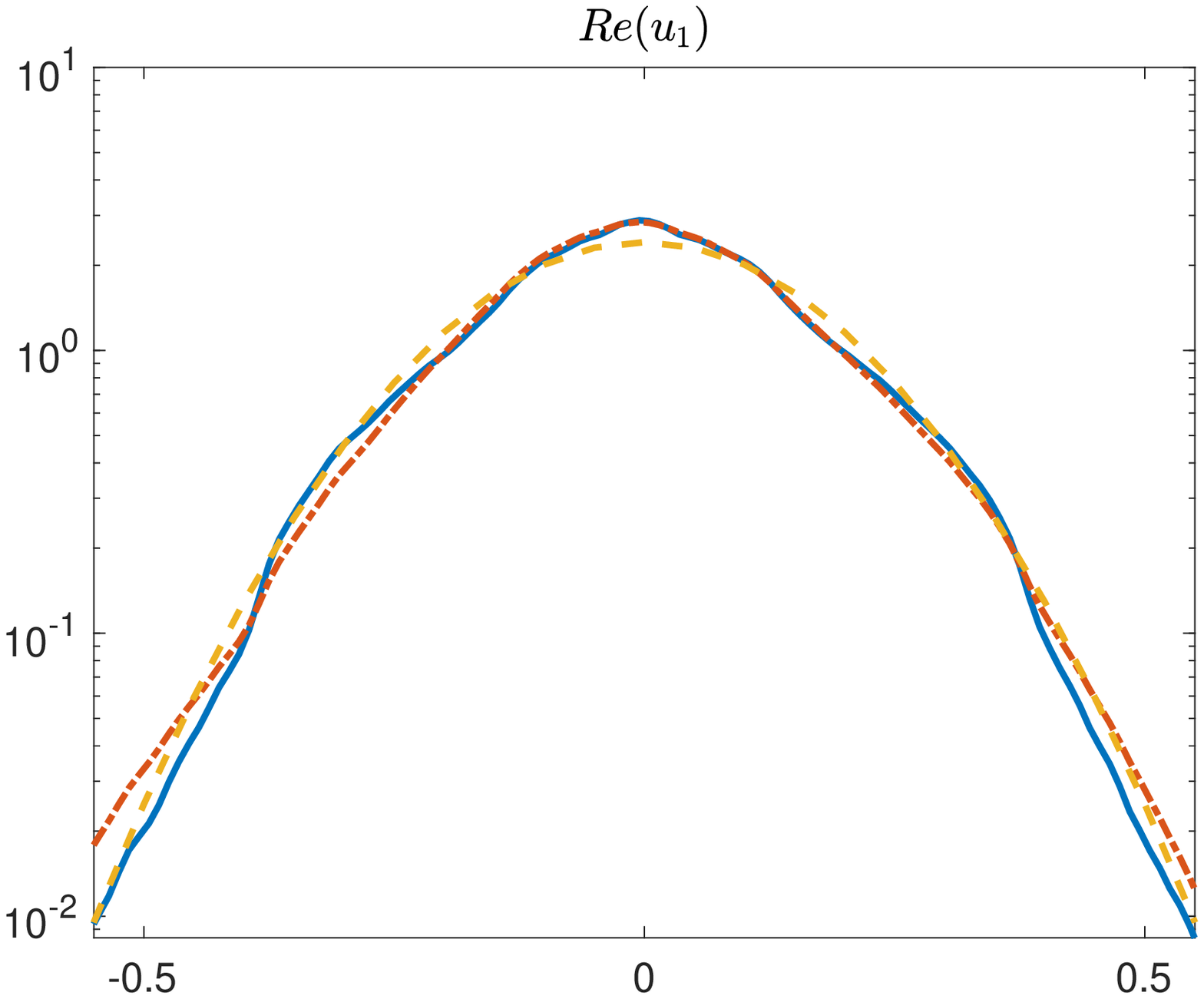}
	\includegraphics[width=.49\textwidth]{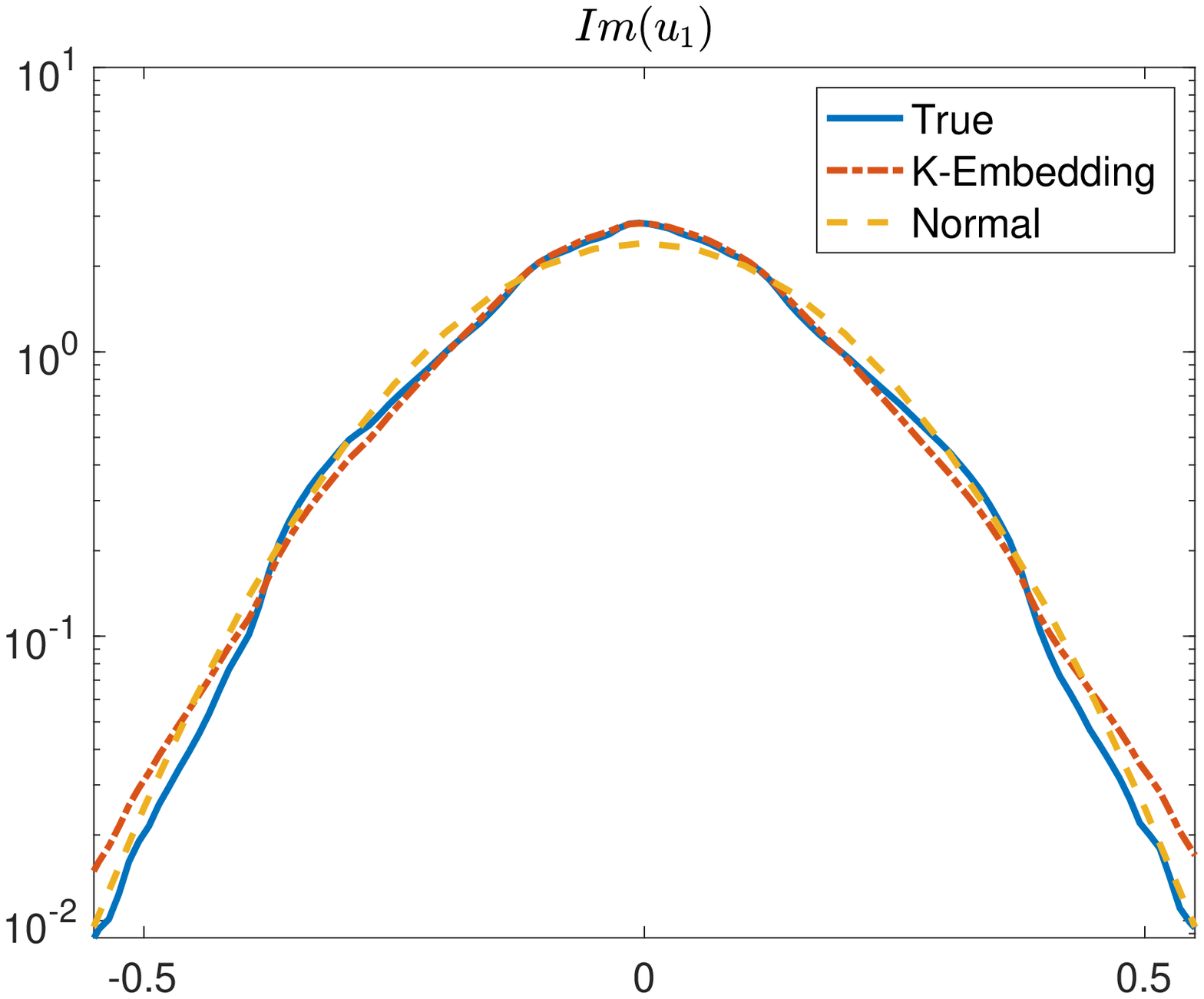}
	\caption{The kernel embedding approximation of the marginal distribution of $x$ in \eqref{eq:real_rep}. The blue solid curves are the PDFs obtained from the observed time series of $u^{N}_{1}$ using kernel density estimates for smoothing. The red dash curves are the order-$5$ kernel embedding estimates. The yellow curves are the density function of the normal distribution of the same mean and variance.}
		\label{fig:KS_2}
	\end{center}
\end{figure}

Figure~\ref{fig:KS_2} shows the performance of the kernel embedding in estimating $\hat{\rho}(x)$. One can see that, compared with the Gaussian approximation, the kernel embedding estimator has a better performance near the peaks of the target density functions. For the corresponding parameter estimation problem, with the Lyapunov equation \eqref{eq:lya_eq}, the parameter space of \eqref{eq:KS_red} can be reduced to
\begin{equation}\label{eq:KS_para}
\theta : = \left(\Lambda, \Gamma \right),
\end{equation}
where $\Lambda$, as a covariance matrix of the auxiliary variable $v$, is SPD and $\Gamma$ is positive definite. Connecting to Section~\ref{sec:scheme}, components of $\Lambda$ are the equilibrium parameters, $\theta_{eq}$. 

\subsection{Parameter estimation and numerical results}

Next, we will estimate the parameters $\theta=(\Lambda,\Gamma)$ using the scheme introduced in Section~\ref{sec:scheme} and check whether the imperfect model in \eqref{eq:KS_red} with the estimated parameters can produce a better estimate of the full response statistics compared to the linear responses as shown in Figure~\ref{fig:KS_1}. First, notice that the imperfect model \eqref{eq:KS_red} yields an equilibrium distribution \eqref{eq:peq_red_KS}, which suggests the auxiliary variable $v$ is Gaussian with covariance matrix $\Lambda$ at the equilibrium. Thus, we can use the sampled covariance of $\dot{x}$, approximated through a finite-difference, to directly estimate $\Lambda$. We will denote the estimate for $\Lambda$ as $\hat{\Lambda}$.

\begin{figure}[ht!]
\begin{center}
\includegraphics[width=.49\textwidth]{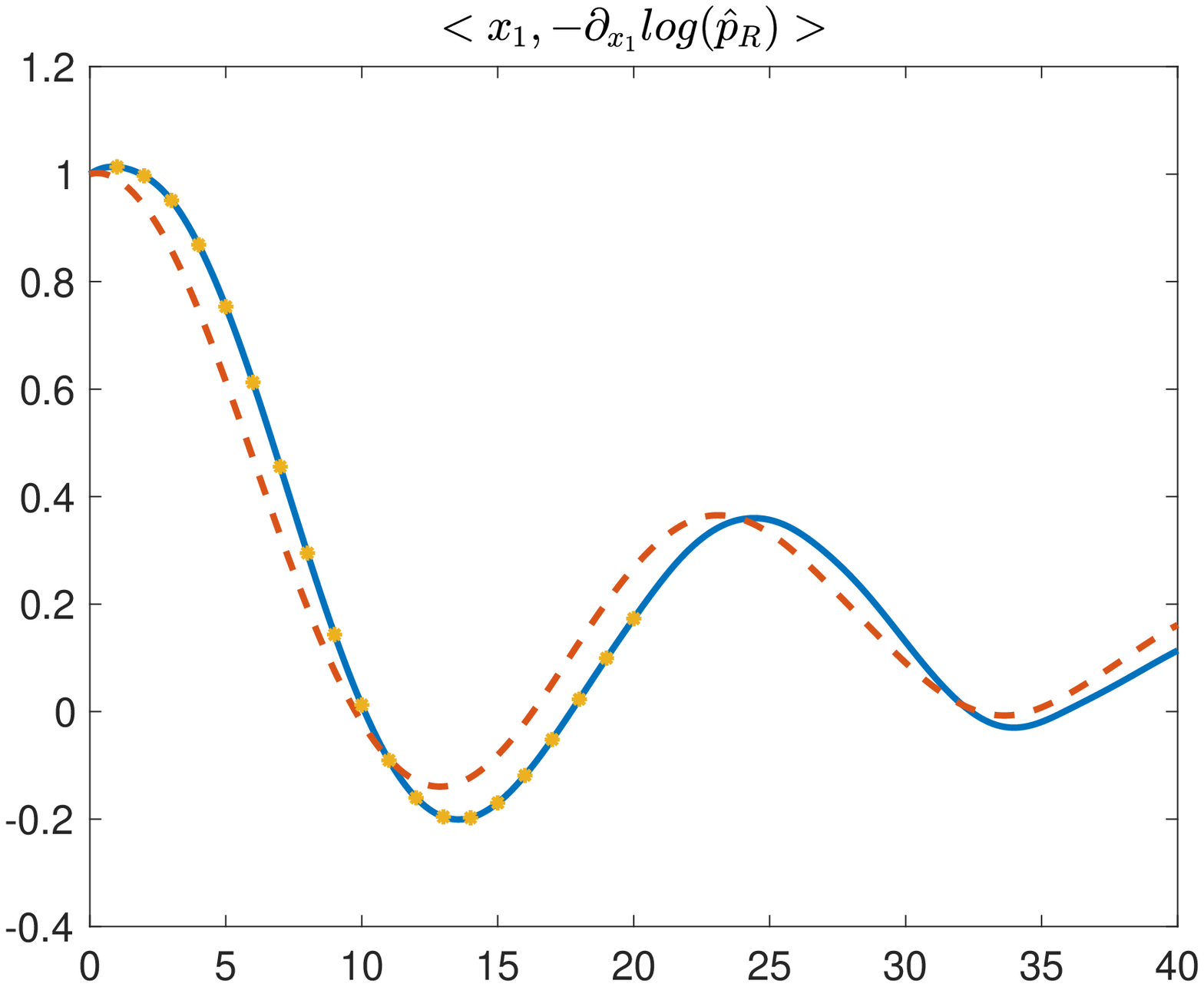}
\includegraphics[width=.49\textwidth]{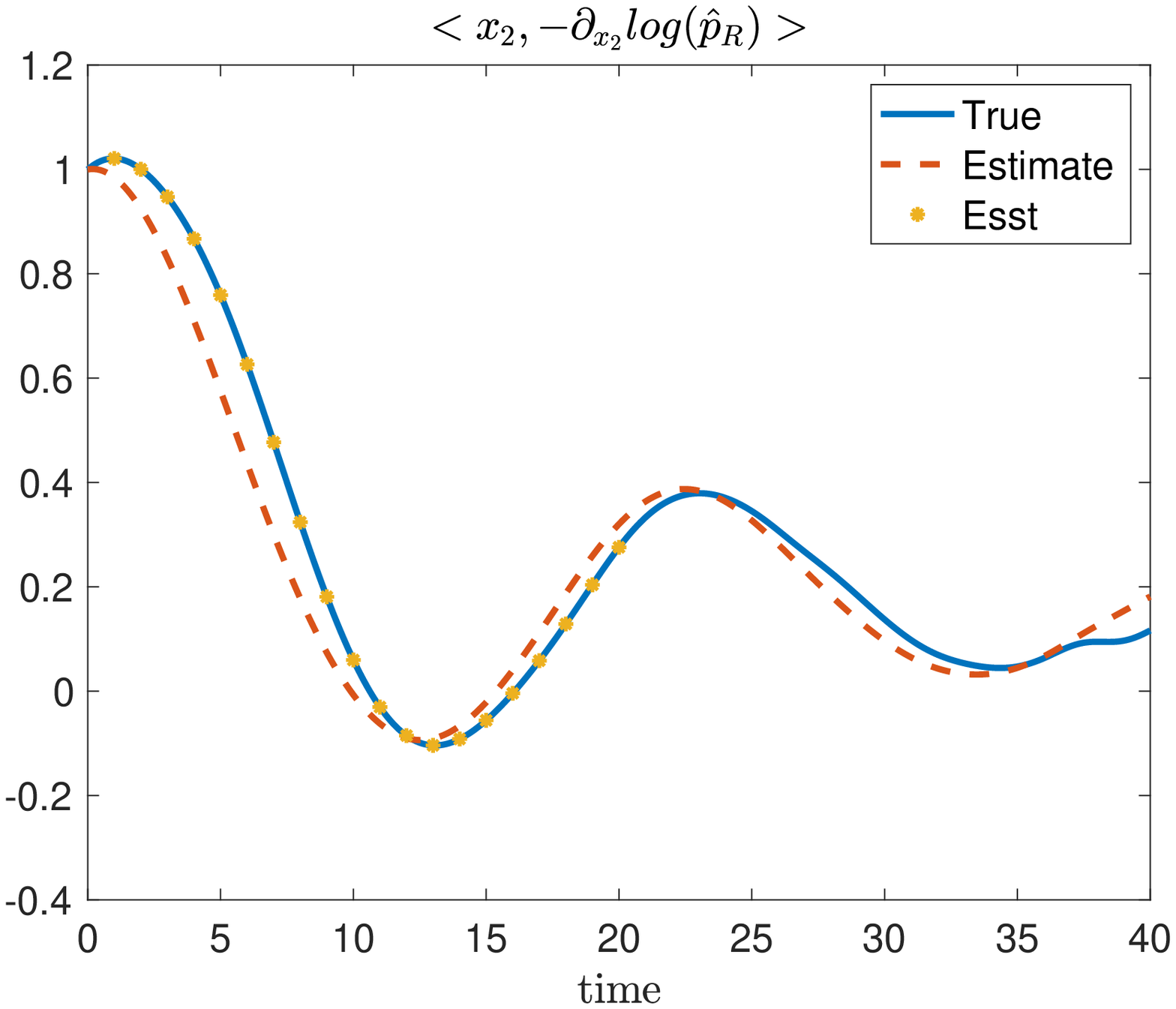}
\caption{Result of the nonlinear least squares fitting \eqref{eq:KS_nls}. The values of the essential statistics $\{g^{\dagger}_{j}(t_{i})\}$, taken from the blue curves, are highlighted by yellow stars. The red dash curves are the corresponding two-time statistics of the imperfect model \eqref{eq:KS_red} the estimated parameters reported in Table~\ref{tab:KS_est}. }
\label{fig:KS_3}
\end{center}
\end{figure}

For the damping coefficient $\Gamma$, we take the diagonal elements of the approximated linear response operator \eqref{eq:lin_resp_KS1} as the essential statistics (they are shown as the yellow stars in Figure~\ref{fig:KS_3}), and we are able to set up the following nonlinear least square problems to estimate $\Gamma$:
\begin{equation} \label{eq:KS_nls}
\min_{\Gamma} \sum_{i = 1}^{K} \sum_{j=1}^2 \left(g^{\dagger}_{j}(t_i) - \hat{g}_{j}(t_{i};\Gamma,\hat{\Lambda})\right)^2,
\end{equation}
where
\begin{equation*}
g^{\dagger}_{j}(t) = \mathbb{E}_{p^{\dagger}_{eq}}\left[x_{j}(t)\partial_{x_{j}} \log(\hat{p}_{R,M}) \right], \quad \hat{g}_{j}(t;\Gamma,\hat\Lambda) = \mathbb{E}_{\rho_{eq}(\Gamma,\hat{\Lambda})}\left[x_{j}(t)\partial_{x_{j}} \log(\hat{\rho}) \right], \quad j = 1,2.
\end{equation*}
In particular we take $t_i = i$ for $i = 1,2,\dots, 20$. Here, the step length of the marginal essential statistics is larger than that of the MD example in \eqref{eq:Lan_nls} since the correlation times of the linear response operator are different, as suggested by Figure~\ref{fig:corr_fit} and Figure~\ref{fig:KS_3}. We should point out that the $\mathbb{E}_{\peq^\dagger}$ will be approximated by Monte-Carlo over the true time series of $\{u_{1,t},u_{2,t},g_t\}_{t=1,\ldots,T}$. We use Euler-Maruyama method to generate time series from the imperfect dynamics at $\Gamma$ and $\hat{\Lambda}$. Similar to the setup of the numerical scheme of the underlying dynamics \eqref{eq:Fourier_full}, we drop the first $5 \times 10^4$ time units to relax the system, and sub-sample $x$ in the next $5 \times 10^5$ time units with step length $\Delta t = 0.1$. Table \ref{tab:KS_est} records the estimated parameter values.

\begin{figure}[ht!]
	\begin{center}
         \includegraphics[width=.49\textwidth]{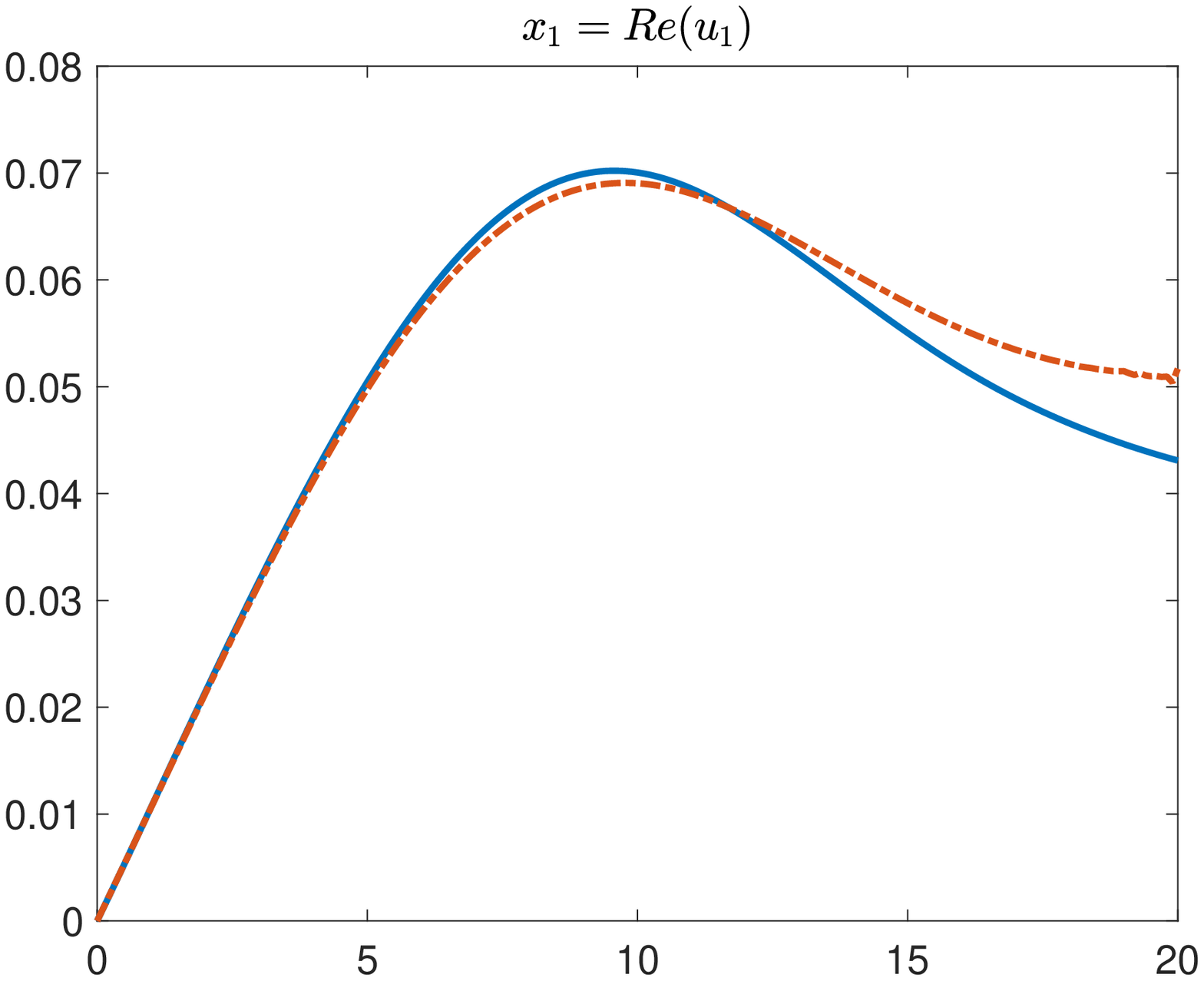}
         \includegraphics[width=.49\textwidth]{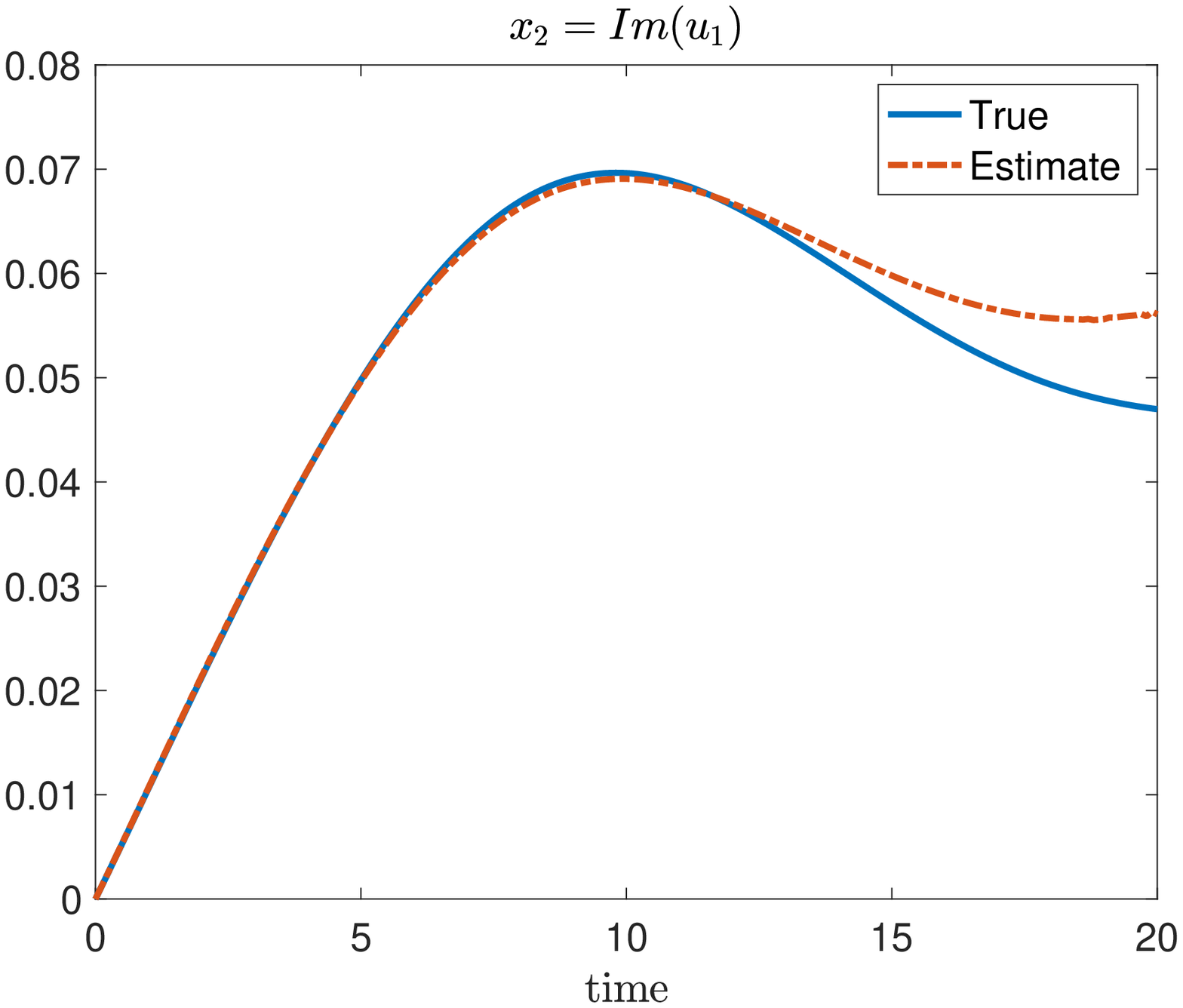}
		\caption{The prediction of the full response. The blue solid curves indicate the same full response in Figure \ref{fig:KS_1}; The red dash curves are the full response  from the imperfect model \eqref{eq:KS_red} based on the parameter estimates reported in Table~\ref{tab:KS_est}.}
		\label{fig:KS_4}
	\end{center}
\end{figure}

\begin{table}[ht]
	\begin{center}
		\begin{tabular}{|c|c|c|c|c|c|c|}\hline
			$ \Lambda_{11} $& $\Lambda_{12}$ & $\Lambda_{22}$ & $\Gamma_{11}$ & $\Gamma_{12}$ & $\Gamma_{21}$ & $\Gamma_{22}$  \\
			\hline
			 $1.3675\times 10^{-3}$& $-1.01 \times 10^{-5}$ &  $1.3544 \times 10^{-3}$ & 0.5327 &  -0.5038 &  -0.4953 & 0.6176\\ \hline  
		\end{tabular}
		\caption{Estimates of the parameter in \eqref{eq:KS_para}.}
		\label{tab:KS_est}
	\end{center}
\end{table}

In Figure \ref{fig:KS_3}, one can see that the estimated linear response statistics (red dashes) from the imperfect model reproduces the qualitative feature of the true linear response statistics. Figure \ref{fig:KS_4} shows the recovery of the full response using the imperfect model in \eqref{eq:KS_red} based on the estimated parameters in Table~\ref{tab:KS_est}. One can see that the real and imaginary components of $u_1$ share similar behavior in both the two-point statistics and the full response statistics. We note that the imperfect model is able to capture such symmetry. Compared with the approximated linear response statistics showed in Figure \ref{fig:KS_1}, the full response statistics of the imperfect model produce a more accurate estimation for larger values of $t$.

\section{Summary and discussion} \label{sec:sum}

Motivated by FDT, we have demonstrated a response-based estimation method to determine parameters that are responsible for the response statistics in stochastic dynamics. Unlike the prior works  \cite{HLZ:17,HLZ:19}, which were focused on the perfect model case, the emphasis of this paper is on a more practical scenario where the model error is present. We started with a conventional two-scale model and showed that the response-based approach can outperform a classical averaging method that is based on the equilibrium statistics, even in the regime where the averaging theory is valid. 

More importantly, we demonstrated that in general, the parameter estimation problems are far more complicated by the fact that only partial observations are available and the response statistics might not be accessible. To overcome these difficult issues, we have outlined a general strategy that can be efficiently implemented in practice. In particular, we introduced the concepts of admissible external forcing, for which the forces on the variable of interest ($X$) can be identified and computed, as well as the marginal linear response (MLR) statistics, as computable statistics in place of the linear response statistics of the underlying full dynamics that are inaccessible in general.

To illustrate this general strategy in a specific context, we considered the Langevin dynamics model, a widely used stochastic model in practice. We studied two examples, including a one-dimensional chain model, coarse-grained from a full MD model; and an extended Langevin dynamics, obtained from a severely truncated KS equation in the Fourier space.  In both cases, we demonstrate how our general strategy is implemented and validated. More specifically,  after the model parameters are determined from the least-squares approach, we compared the {\it nonlinear} responses of the imperfect model with those of the underlying dynamics under admissible forcing. A quantitative agreement has been found. 

For the second example, we also employed a nonparametric kernel embedding formulation to estimate the marginal density, and more importantly, the marginal essential statistics, to enable the estimation using the response statistics. We report the convergence analysis of this estimator in \cite{zhl:19b}.

\section*{Acknowledgment}

This research was supported by NSF under grant DMS-1619661 and DMS-1819011. JH also acknowledges the support from NSF grant DMS-1854299.

\appendix

\section{Approximation of the full response statistics via Monte-Carlo} \label{app:A}

To test the validity of the parameter estimates, we compare the full response statistics of the underlying dynamics and those of the imperfect model. For example, in Section \ref{sec:chain}, to generate Figure \ref{fig:Lan_resp}, one needs to compute the full response statistics of the underlying dynamics \eqref{eq:Lan_temp} and the imperfect dynamics \eqref{eq:CG_sgd} under a constant external forcing on the first block of particles and the first CG particle, respectively. In this Appendix, we review the Monte-Carlo (MC) approach used in approximating the full response statistics $\Delta \mathbb{E}[A](t)$ \eqref{eq:full_resp} under a general It\^o diffusion setting.

Recall that in Section \ref{sec:review}, we have used \eqref{eq:Ito} and \eqref{eq:Ito_per} to represent the unperturbed and perturbed dynamics respectively. In the verification, the unknown parameter $\theta$ in the two dynamics will be replaced by the estimate $\hat\theta$. To compute $\Delta \mathbb{E}[A](t)$, since $\mathbb{E}_{p_{eq}}[A(X)]$ in \eqref{eq:full_resp} is an equilibrium statistics of the unperturbed dynamics, it can be directly approximated by an ensemble average
\begin{equation}\label{eq:eq_stat}
    \mathbb{E}_{p_{eq}}[A(X)] \approx \frac{1}{N} \sum_{i =1}^{N} A(x_i),
\end{equation}
where $\{x_{i} = X(t_{i})\}$ denotes a time series of the unperturbed dynamics \eqref{eq:Ito} generated at the equilibrium state. Thus, we are able to reduce our problem to approximating the response statistics $\mathbb{E}_{p^{\delta}}[A(X^{\delta}(t))]$.

Introducing the transition kernel $\rho^{\delta}(x,t|x_0,0)$ of the perturbed dynamics \eqref{eq:Ito_per}, $\mathbb{E}_{p^{\delta}}[A(X^{\delta}(t))]$ can be calculated by a double integral
\begin{equation}\label{eq:dou_int}
    \mathbb{E}_{p^{\delta}}[A(X^{\delta}(t))] = \int \int A(x) \rho^{\delta}(x,t|x_0,0) p_{eq}(x_0) \td x \td x_0,
\end{equation}
where we have used the fact that the perturbed dynamics \eqref{eq:Ito_per} is also initiated at the equilibrium state of the unperturbed dynamics \eqref{eq:Ito} with equilibrium distribution $p_{eq}$. Further notice that the transition kernel satisfies the Fokker-Planck equation
\begin{equation*}
    \frac{\partial \rho^{\delta}}{\partial t} = (\mathcal{L}^{\delta})^{*}\rho^{\delta}, \quad \rho^{\delta}(x, 0 |x_{0},0) = \delta(x-x_0),
\end{equation*}
where $\mathcal{L}^{\delta}$ denotes the generator of the perturbed dynamics \eqref{eq:Ito_per}. One can rewrite the transition kernel using the semi-group notation: $\rho^{\delta}(x,t|x_0,0) = e^{t(\mathcal{L}^{\delta})^{*}}\delta(x-x_0)$.  Then, the double integral \eqref{eq:dou_int} becomes
\begin{equation}\label{eq:MC_1}
    \begin{split}
    \mathbb{E}_{p^{\delta}}[A(X^{\delta}(t))]  & = \int \int A(x) e^{t(\mathcal{L}^{\delta})^{*}}\delta(x-x_0) p_{eq}(x_0) \td x \td x_0 \\
    & = \int e^{t\mathcal{L}^{\delta}} A(x) \int \delta(x-x_0) p_{eq}(x_0) \td x_0 \td x \\
    & = \int e^{t\mathcal{L}^{\delta}} A(x)p_{eq}(x)\td x.
     \end{split}
     \end{equation}
Here, $u^{\delta}(x,t): = e^{t\mathcal{L}^{\delta}}A(x)$ solves the Kolmogorov's backward equation (dual form of the Fokker-Planck equation)
\begin{equation*}
    \frac{\partial u^{\delta}}{\partial t} = \mathcal{L}^{\delta}u^{\delta}, \quad u^{\delta}(x,0) = A(x),
\end{equation*}
which suggests that $u^{\delta}(x,t)$ corresponds to the conditional expectation $ \mathbb{E}\left[A(X^{\delta}(t)) \big| X^{\delta}(0) = x \right]$. Together with Eq. \eqref{eq:MC_1}, we derive the following two-layer Monte-Carlo approximation,
\begin{equation}\label{eq:MC}
    \mathbb{E}_{p^{\delta}}[A(X^{\delta}(t))] \approx \frac{1}{N} \sum_{i=1}^{N} u^{\delta}(x_i,t) \approx \frac{1}{MN}\sum_{i=1}^{N}\sum_{j=1}^{M} X^{\delta,x_{i}}(t; U_{j}),
\end{equation}
where $X^{\delta,x}(t; U_{j})$ denotes the realization of \eqref{eq:Ito_per} with initial condition $X^{\delta}(0)=x_{i}$, given a sample path $U_{j}$ of the Wiener process $U_{s}$ ($0\leq s \leq t$) in \eqref{eq:Ito_per}.

In the numerical implementation, after generating the time series $\{x_i\}$ at $p_{eq}$, one can duplicate $x_i$ up to $M$ times and solve \eqref{eq:Ito_per} in parallel based on the samples from the time series as initial conditions. In particular, we take $N = 10^{5}$ and $M=10^{3}$ in estimating the full response statistics of the underlying dynamics and imperfect model shown in Figure \ref{fig:Lan_resp}.

\section{Heat flux and the density ansatz} \label{app:ansatz}

The main goal of this Appendix is to deduce a parametric form of the equilibrium density \eqref{eq:CG_rho1} of the coarse-grained system. A system with non-uniform temperatures usually exhibits heat conduction, in which case the induced heat flux becomes relevant in our parameter estimation procedure (e.g., \cite{lepri2003thermal,wu2015consistent}).
In Section~\ref{sec:chain}, the local heat flux $j_i$ corresponding to the definition of the CG variables  $(\bm{\xi}, \bm{p})$ (\ref{eq:CG_par},\ref{eq:CG_obs}) is given by,
\begin{equation*}
    j_i \propto (p_{i+1} + p_{i}) V_{i}'(\xi_i),
\end{equation*}
where $V_i$ denotes the local potential (e.g., \cite{lepri2003thermal}).

Let $\rho$ be the equilibrium distribution of the CG variables in \eqref{eq:CG_sgd}. It is important to note that this notion of equilibrium density is more general due to the non-zero current. Such a density function $\rho$ corresponds to the steady-state solution of the Fokker-Planck equation of the CG model. Except for very special cases, the explicit form of $\rho$ is unknown. Due to the presence of heat conduction, it is necessary to incorporate the heat flux in approximating the density function $\rho$.

Motivated by the maximum entropy approach, we consider the following two set of constraints on the local Hamiltonian $\mathcal{H}_i:= \frac{p_{i}^2}{2} + V_i$ and $j_i$,
\begin{equation}\label{eq:const_1}
    \mathbb{E}_{\rho}[\mathcal{H}_i] = \mathbb{E}_{p_{eq}^{\dagger}}[\mathcal{H}_i], \quad \mathbb{E}_{\rho}[j_i] = \mathbb{E}_{p_{eq}^{\dagger}}[ j_i], \quad i = 1,2,\dots, n_b,
\end{equation}
which can be interpreted as a system of integral equations to $\rho$. Each heat flux constraint in \eqref{eq:const_1} involves a three-dimensional integral with respect to $\xi_i$, $p_{i}$, and $p_{i+1}$, which causes issues in our later applications of maximum entropy. We resolve such an issue by exploiting the equality using the stationarity,
\begin{equation*}
    \frac{\td }{\td t} \mathbb{E}_{p_{eq}^{\dagger}}[V_{i}(\xi_i)] = 0, \quad i = 1,2,\dots, n_b,
\end{equation*}
which implies that
\begin{equation*}
    \mathbb{E}_{p_{eq}^{\dagger}}\left[p_i V_{i}'(\xi_i)\right] =  \mathbb{E}_{p_{eq}^{\dagger}}\left[p_{i+1} V_{i}'(\xi_{i})\right].
\end{equation*}
As a result, the constraints in \eqref{eq:const_1} can be equivalently rewritten as,
\begin{equation}\label{eq:const_2}
    \mathbb{E}_{\rho}[\mathcal{H}_i] = \mathbb{E}_{p_{eq}^{\dagger}}[\mathcal{H}_i], \quad \mathbb{E}_{\rho}[p_{i}V_{i}'(\xi_{i})] = \mathbb{E}_{p_{eq}^{\dagger}}[p_{i}V_{i}'(\xi_{i})], \quad i = 1,2,\dots, n_b,
\end{equation}
such that each constraint is a two-dimensional integral equation. With constraints in \eqref{eq:const_2}, the corresponding maximum entropy solution is of the form
\begin{equation}\label{eq:CG_rho2}
    \rho(\bm{\xi}, \bm{p}) \propto \exp\left[ - \sum_{i=1}^{n_b} b_i \left(\frac{p_{i}^2}{2} + V_{i}(\xi_i)\right) - \sum_{i=1}^{n_b} C_i p_{i} V_{i}'(\xi_i) \right],
\end{equation}
together with the periodic boundary condition. Such an ansatz in \eqref{eq:CG_rho2} can be interpreted as a perturbed Gibbs measure \cite{zwanzig2001nonequilibrium}.

The potential function $V_{i}$ in $\rho$ \eqref{eq:CG_rho2} is assumed to be a quartic function of the following form,
\begin{equation}
    V_{i}(\xi_i) = \sum_{j=1}^4 A_{j,i}\xi_i^{j}, \label{qpotential}
\end{equation}
with unknown parameters $\{A_{j,i}\}$. Substituting \eqref{qpotential} to \eqref{eq:CG_rho2}, we obtain \eqref{eq:CG_rho1} with parameters: $a_{j,i} = A_{j,i}b_i$ for $j=1,\ldots, 4$, and $c_{j,i} = (j+1)A_{j+1,i}C_i$ for $j=1,\ldots,3$. The restriction of parameters makes the maximum entropy formulation to be inconvenient in finding $b_i$ and $c_i$. To overcome this issue, we keep the same ansatz for the equilibrium density in \eqref{eq:CG_rho1} and solve the maximum entropy problem in \eqref{eq:est_a}-\eqref{eq:mom_con} for each coarse grained component.

\section{The isothermal one-dimensional molecular model} \label{app:iso}

In Section~\ref{sec:chain}, we emphasized on the parameter estimation problem of the imperfect model of a coarse-grained molecular model \eqref{eq:Lan_temp} with non-uniform temperature. We also tested the isothermal case, where $k_{B}T_i \equiv k_{B}T = 0.25$ in \eqref{eq:Lan_temp} for $i = 1,2,\dots, N$. We included the results in this section for interested readers. 
As a simplification of \eqref{eq:Lan_temp}, $(\mbr, \mbv)$ obeys the Langevin equation of motion,
\begin{equation}\label{eq:Lan_full}
\begin{split}
\dot{\mbr} &= \mbv,\\
\dot{\mbv} &= - \nabla_{\mbr}U(\mbd) - \gamma \mbv +\sqrt{2k_{B}T\gamma}\dot{W}_{t},
\end{split}
\end{equation}
The Langevin dynamics is a mean to sample canonical ensemble \cite{leimkuhler2016molecular}, or the Gibbs measure. The equilibrium distribution of the relative displacement $\mbd$ and the velocity $\mbv$, is of the form
\begin{equation}\label{eq:p_eq}
p_{eq}(\mbd, \mbv) \propto \exp\left[-\frac{1}{k_{B}T}\left(U(\mbd)+\frac{1}{2}\mbv^{2}  \right)\right],
\end{equation}
under the constraint induced by the periodic boundary condition.

In this isothermal case, the CG model in \eqref{eq:CG_sgd} reduces to
\begin{equation}\label{eq:Lan_red}
\begin{split}
\dot{\mbq} &= \mbp,  \\
\dot{\mbp} &= - \nabla_{\mbq} V(\bm{\xi}; \bm{a}) - \Gamma \mbp + \sqrt{2\beta^{-1}}\sigma \dot{U}_t,
\end{split} 
\end{equation}
where the damping matrix $\Gamma$ satisfies $ \sigma\sigma^{\top} = \Gamma$, and the PMF $V(\bm{\xi}; \bm{a})$, as a function of the CG relatively displacement $\bm{\xi}$ is of the form
\begin{equation}\label{eq:CG_potential}
V(\bm{\xi}, \bm{a}) = \sum_{i = 1}^{J}\sum_{j=1}^4 a_{j}\xi_{i}^{j}, \quad q_{k+n_{b}} = q_{k}, \quad a_{4}>0, 
\end{equation}
with unknown parameters $\bm{a}= (a_{1},a_{2},a_{3},a_{4})$. Adopting the form of the damping matrix $\Gamma$ in \eqref{eq:damping}, the parameter set of the imperfect model \eqref{eq:Lan_red} is
\begin{equation*}
    \theta = (\bm{a}, \beta, \gamma_0, \gamma_1, \gamma_2),
\end{equation*}
with $\theta_{eq}:=(\bm{a}$, $\beta$). To determine $\theta_{eq}$, consider the Gibbs measure of the imperfect model \eqref{eq:Lan_red} (following the notation in \eqref{eq:lin_oper_red}),
\begin{equation}\label{eq:CGp_eq}
\rho_{eq}(\bm{\xi}, \mbp; \theta_{eq}) \propto \exp\left[-\beta\left(V(\bm{\xi};\bm{a})+\frac{1}{2}\mbp^{2} \right)\right],
\end{equation}
which suggests that the temperature $\beta^{-1}$ can be directly estimated from the sample variance of the CG velocities. In practice, since $\{p_{i}\}$ are independent identical distributed at the equilibrium, we take
\begin{equation}\label{eq:est_beta}
\hat \beta^{-1} =  \mathbb{E}_{p_{eq}^{\dagger}}\left[p_{1}^2\right],
\end{equation} 
as the estimate for $\beta^{-1}$. For the parameters in the potential \eqref{eq:CG_potential}, we apply the maximum entropy method to the marginal distribution of $q_{i+1} - q_i$ with respect to $\rho_{eq}$ \eqref{eq:CGp_eq},
\begin{equation*}
\rho(x;\bm{a}, \beta) \propto \exp \left(-\beta \sum\limits_{j = 1}^{4}a_{j}x^{j} \right),
\end{equation*}
 for $i = 1,2,\dots, n_{b}-1$. By maximizing the Shannon's entropy under the moment constraints, we formulate the following maximum entropy estimates for $\bm{a}$,
\begin{equation*}
    \begin{split}
    \hat{\bm{a}} &= \argmin \int -\rho\left(x;\bm{a}, \hat\beta\right) \log \rho\left(x;\bm{a}, \hat\beta\right) \td x, \nonumber \\
\text{s.t. } &\int x^{j}\rho\left(x;\bm{a}, \hat\beta\right) \td x = \frac{1}{n_{b}-1} \sum\limits_{i=1}^{n_{b}-1}\mathbb{E}_{p^{\dagger}_{eq}}\left[\xi_{i}^{j}\right], \quad j = 1,2,3,4,    
    \end{split}
\end{equation*}
where $\hat\beta$ is given by \eqref{eq:est_beta}.

\begin{figure}[ht!]
	\begin{center}
		\includegraphics[width=.49\textwidth]{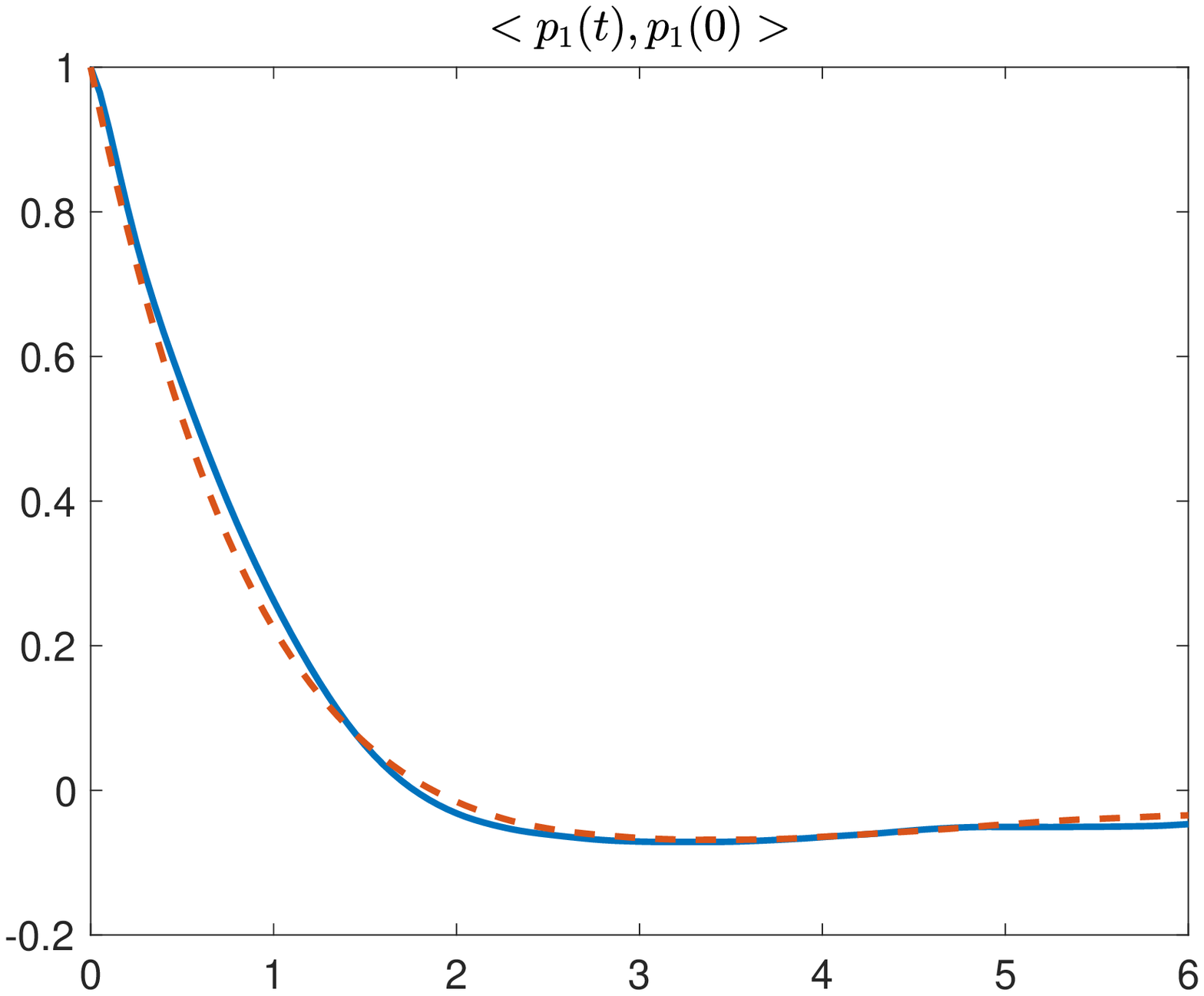}
		\includegraphics[width=.49\textwidth]{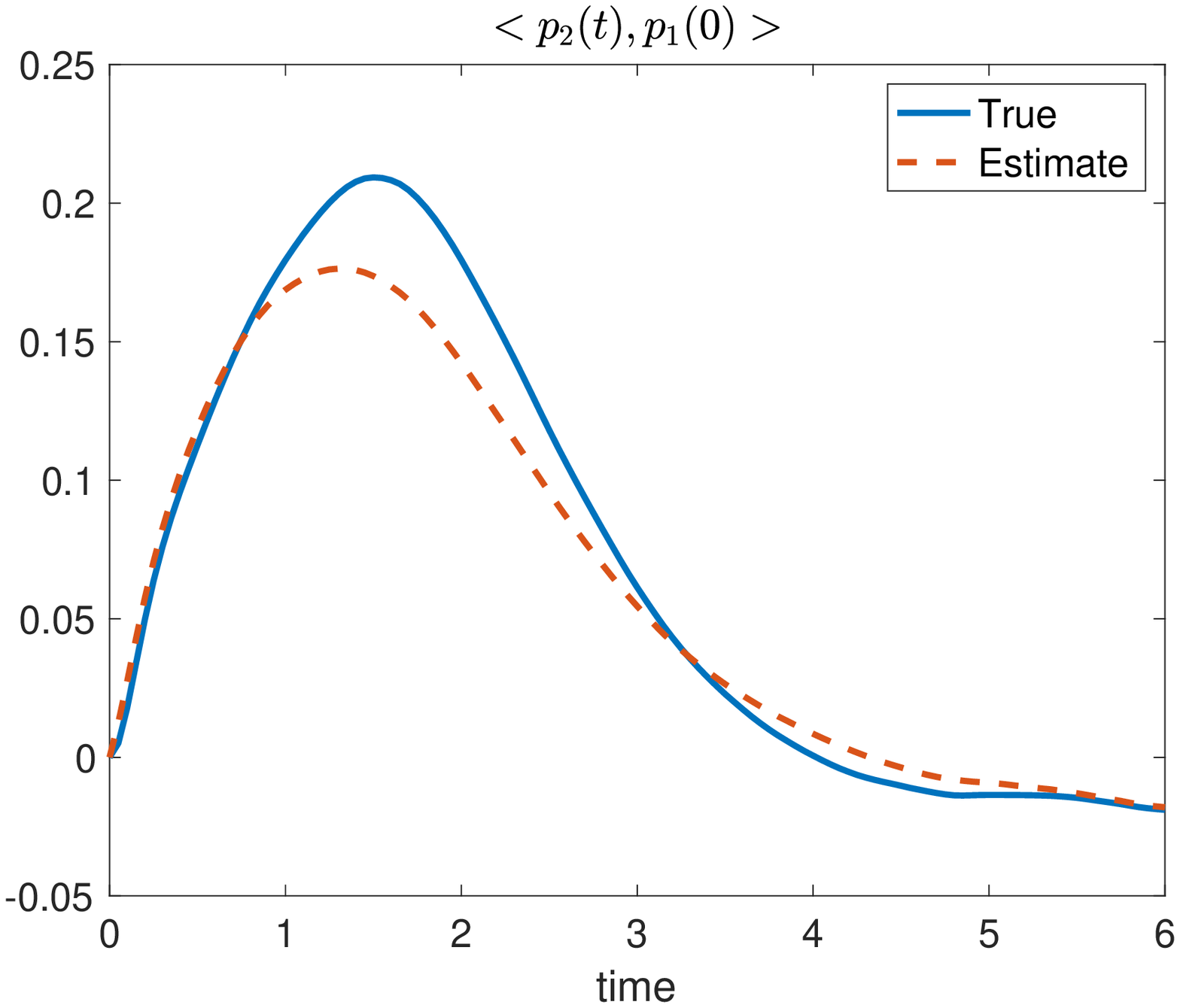}
		\caption{Comparison of the time correlation functions involved in the nonlinear least-squares problem. The blue solids curves are computed from the available CG observations, while the red dash curves are the approximation of the imperfect model based on the estimates.}
		\label{fig:corr_fit_iso}
	\end{center}
\end{figure}

\begin{figure}[ht!]
	\begin{center}
		\includegraphics[width=.49\textwidth]{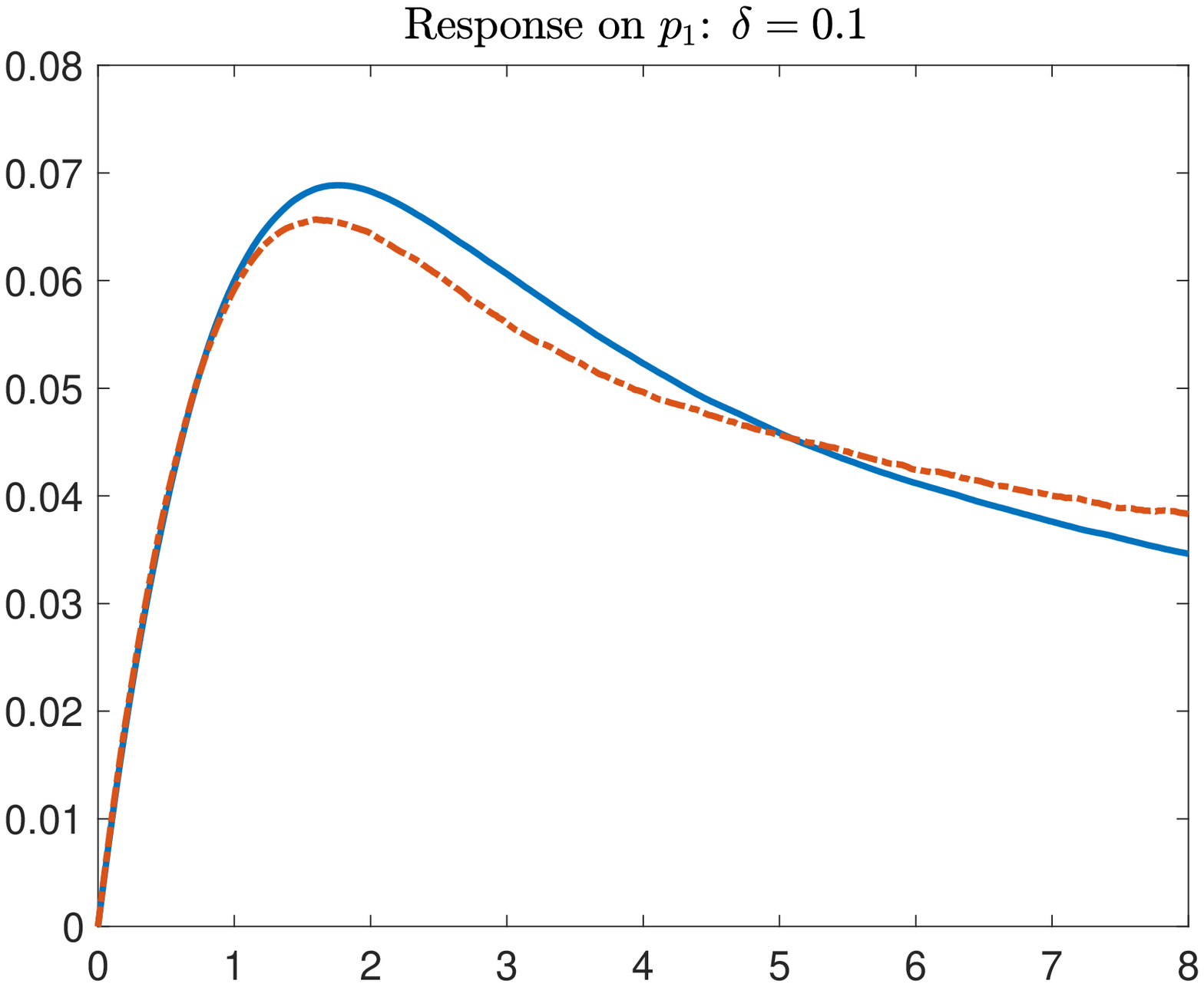}
		\includegraphics[width=.49\textwidth]{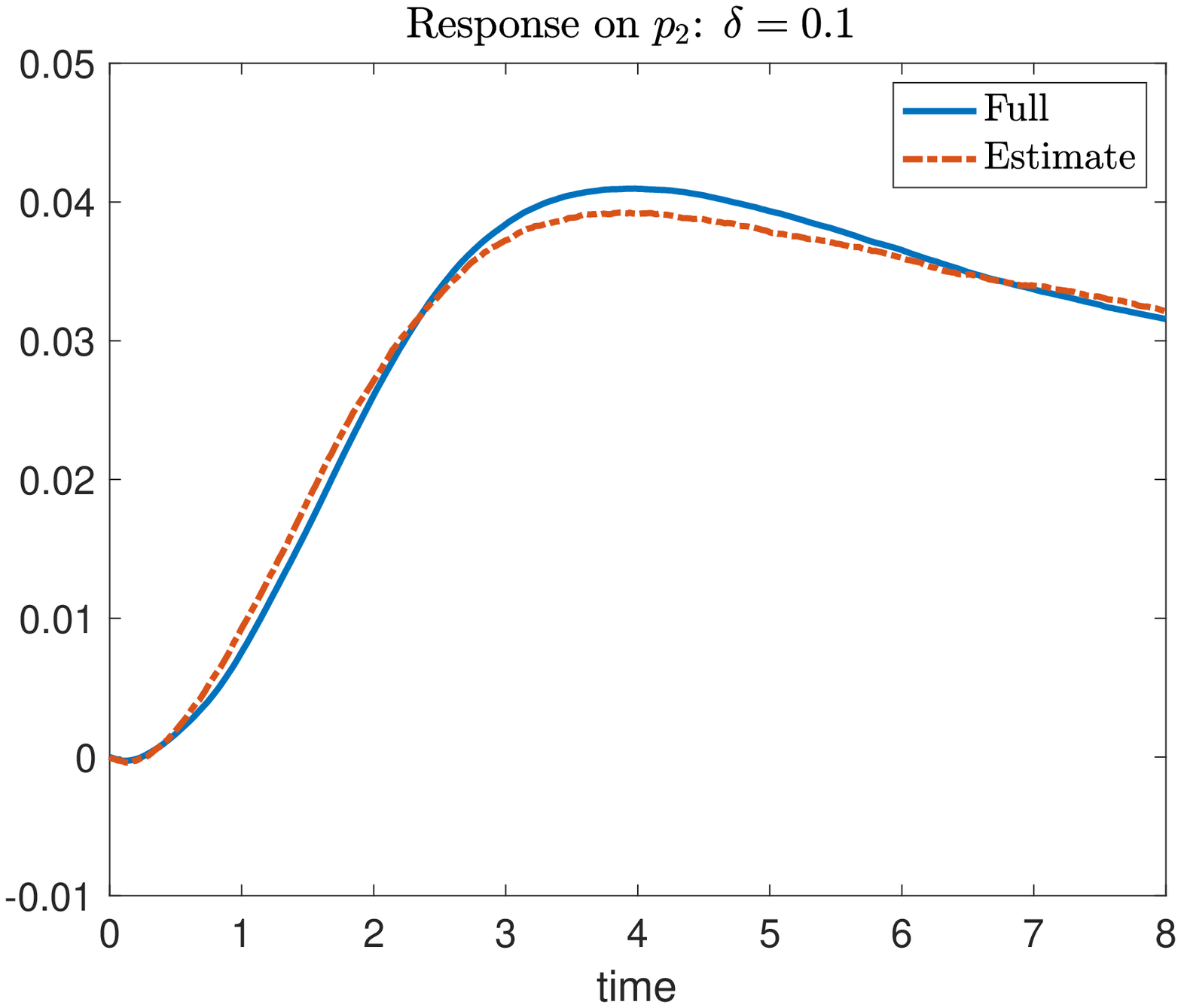}
		\caption{Reconstruction of full response statistics: underlying dynamics (blue solid) v.s. the imperfect model (red dash).}
		\label{fig:Lan_resp_iso}
	\end{center}
\end{figure}

For parameters $\{\gamma_i\}$, we take the same approach as in Section~\ref{sec:chain}. Table \ref{tab:Lan_est_iso} shows the value of the estimates. Figure~\ref{fig:corr_fit_iso} and Figure~\ref{fig:Lan_resp_iso} provide the results in fitting the essential statistics and recovering the full response of the underlying dynamics, respectively.

\begin{table}[ht]
	\begin{center}
		\begin{tabular}{|c|c|c|c|c|c|c|c|c|}
		\hline
			Parameter ($\theta$) & $\beta^{-1}$& $a_{1}$ & $a_{2}$ & $a_{3}$ & $a_{4}$ & $\gamma_{0}$ & $\gamma_{1}$ & $\gamma_{2}$ \\
			\hline
			Estimates ($\hat{\theta})$ & 0.0247 & 0.0193 &  0.148 & -0.0830 &  0.0331 & 1.23 & -0.371 & -0.0420  \\ \hline
		\end{tabular}
	\end{center}
	\caption{Estimates of the parameter in \eqref{eq:Lan_red}. The scale of $a_j$ is different from the scale of $a_{j,i}$ in Table~\ref{tab:Lan_est} since $\beta$ is a prefactor of potential function $V$ in \eqref{eq:CGp_eq}.}
	\label{tab:Lan_est_iso}
\end{table}

As in Section~\ref{sec:chain}, we consider the damping matrix $\Gamma$ in the form of symmetric and circulant matrix with non-zero elements $(\gamma_0, \gamma_1, \gamma_2)$. In our numerical experiments of the isothermal molecular model, we also tried damping matrix $\Gamma$ with smaller bandwidth. We found that using $\Gamma = \gamma_0 I$ produces a much less accurate recovery of the time correlation (results are not reported) compared to those obtained using non-diagonal $\Gamma$ as presented in Figure~\ref{fig:corr_fit_iso}.  We also did not continue with testing the damping matrix with larger bandwidth since the estimates of $\gamma_2$ are already small.

\bibliographystyle{unsrt}

\bibliography{aipsamp,cgmd}

\end{document}